\documentclass[a4paper,10pt]{article}
\usepackage[left=3cm,right=3cm,top=3cm,bottom=3.5cm]{geometry}
\usepackage[latin1]{inputenc}
\usepackage{amscd}
 \usepackage{amsthm}
\usepackage{amsfonts}
\usepackage{color}
\usepackage{multicol}
\usepackage{enumerate}
\usepackage{makeidx}
\usepackage{indentfirst}
\usepackage{comment}
\usepackage{amsmath,amssymb}
\usepackage{authblk}

\newcommand{\fz}{\mathfrak{z}}
\newcommand{\zw}{b}   
\newtheorem{theorem}{Theorem}[section]
\newtheorem{proposition}[theorem]{Proposition}
\newtheorem{corollary}[theorem]{Corollary}

\newtheorem{lemma}[theorem]{Lemma}
\theoremstyle{definition}
\newtheorem{remark}{Remark}[section]
\newtheorem{definition}{Definition}[section]
\newtheorem{counter-example}[theorem]{Counter-example}
\newcommand{\com}[1]{%
  \ifnum\visibleComment=0
  \else {\bf #1}
  \fi
}
\newcommand{\visibleComment}{0}
\title{No infinite spin for partial collisions converging to isolated central configurations on the plane}

\author[1]{Anna Gierzkiewicz\thanks{anna.gierzkiewicz@uj.edu.pl}}
\author[2]{Rodrigo G. Schaefer\thanks{rodrigo.schaefer@enti.cat}}
\author[1]{Piotr Zgliczy\'nski\thanks{umzglicz@cyf-kr.edu.pl}}
\affil[1]{Jagiellonian University,
ul. \L ojasiewicza 6, 30-348 Krak\'ow, Poland}
\affil[2]{Escola de Noves Tecnologies Interactives (ENTI-UB), Carrer de Baldiri Reixac, 7, 08028, Barcelona}
\date{%
    \today
}

\begin{document}

\maketitle

\begin{abstract}
In the $n$-body problem, when a~cluster of bodies tends to a collision, then its normalized shape curve converges to the set of normalized central configurations, which has $SO(2)$
symmetry in the planar case. This leaves a possibility that the normalized  shape curve  tends to the circle obtained by rotation of some central configuration instead of a particular point on it. This is the \emph{infinite spin problem} which concerns the rotational behavior of total collision orbits in the $n$-body problem. The question also makes sense for  partial collision.

We show that the infinite spin is not possible if the limiting circle is isolated from other connected components of the set of normalized central configurations. Our approach extends the method from recent work for total collision by Moeckel and Montgomery, which was based on a combination of the center manifold theorem with {\L}ojasiewicz inequality. To that we add a shadowing result for pseudo-orbits near normally hyperbolic manifold and careful estimates on the influence of other bodies
on the cluster of colliding bodies.
\\[2mm]
\\[2mm]
\end{abstract}

\section{Introduction}

We study the planar $n$-body problem, that is, the description of the motion of $n$ point particles with masses $m= (m_1, \dots,m_n)$ and positions $q = (q_1, \dots, q_n)\in \mathbb{R}^{2n}$ on the plane under their mutual Newtonian gravitational forces modeled by the equations
\begin{equation}\label{eq:second-order-formulation}
m_i \ddot{q}_i = \frac{\partial U}{\partial q_i}(q)
\text{, \quad where \quad }
U(q):=\sum_{i<j}\frac{m_i m_j}{q_{ij}}
\text{ \quad and \quad }
q_{ij}:=|q_i-q_j|.
\end{equation}
We are particularly interested in a situation when a
\emph{singularity collision} occurs, when there exists some finite time $T\in\mathbb{R}$ such that if we denote
\begin{equation*}
L = (L_1, \dots, L_n) := \lim_{t\to T} q(t)\text{,}
\end{equation*}
then $L\in \Delta$, where $\Delta$ is the \emph{collision set}, namely
\begin{equation*}
    \Delta:= \bigcup\limits_{1\leq i< j \leq n} 
    \{L=(L_1, \dots, L_n)\in \mathbb{R}^{2n} :\: L_i=L_j\}.
\end{equation*}

First, we separate the (indexes of) bodies into \emph{clusters}:
the cluster $G$ is defined to be $G := \left\{i: q_i\rightarrow L_{G} \text{ as }t\rightarrow T\right\}$.
A collision occurs if for some cluster has cardinality $\#G  = k >1$. We have a \emph{total collision} if the colliding cluster contains all $\#G  = n$ bodies. Otherwise, a collision is called a \emph{partial collision}.
We will also say that two bodies $q_i$ and $q_j$ \emph{collide at $T$} if their limit sets coincide: $L_i = L_j$.

Let us single out one cluster $\mathcal{G}$ with  $1<\# \mathcal{G} =k < n$ and focus on \emph{collision orbits} $q_\mathcal{C}=\{q_i,\: i\in \mathcal{G}\}$.
Let us denote the common limit set $L_\mathcal{G}=L_i$, where $i \in \mathcal{G}$.  Let us denote the center of mass of the cluster $c_\mathcal{G}$ by
\begin{equation}
c_\mathcal{G}=\left( \sum_{i\in \mathcal{G}} m_i \right)^{-1}  \sum_{i\in \mathcal{G}} m_iq_i. \label{eq:cofmassG}
\end{equation}
 Then the quantity
\begin{equation}
  I^0_\mathcal{G}(q)=\sum_{i \in \mathcal{G}} m_i |q_i - c_\mathcal{G}|^2,   \label{eq:IG-def}
\end{equation}
is the \emph{intrinsic moment of inertia of the colliding cluster}. As suggested by the authors in \cite{Moeckel2023}, the quantity $r=\sqrt{I^0_\mathcal{G}(q)}$ is a convenient measure of the distance to a~collision.
The vector
\begin{equation}
\left(\hat{q}_\mathcal{G}(t)\right)_{i \in \mathcal{G}}=\left(\frac{q_i(t)-c_\mathcal{G}(t)}{\sqrt{I^0_\mathcal{G}(q(t))}}\right)_{i \in \mathcal{G}}  \label{eq:qhatC}
\end{equation}
is the corresponding \emph{normalized configuration}. In the case of total collision we will omit subscript $\mathcal{G}$.
We will denote the set of normalized configurations with $k$-bodies as $S_k$. It is immediate that $S_k$ is a sphere in $\mathbb{R}^{2 k}$.

A classical result of Chazy \cite{Ch18} about total collision solutions asserts that their normalized configuration curve $\hat{q}(t)$ converges to the set of normalized central
configurations as $t \to T$. In the case of partial collisions, an analogous result is stated in \cite{Saari71} (see also \cite{Saari84}).

Let us denote by $I(q)$ a quadratic form which is  the moment of inertia of configuration $q \in (\mathbb{R}^d)^n$ with respect to origin by
\begin{equation}
 I(q)=\sum_{i=1}^n m_i |q_i|^2.  \label{eq:I-def}
\end{equation}
As in \cite{Moeckel2023} $I(q)$ will play a role of a square of mass norm used often in the sequel.

\begin{definition}
A point in $q \in \left(\mathbb{R}^2\right)^{n}$  is a \emph{central configuration} (shortly: \emph{CC}) if there exists a~constant $\lambda > 0$ such that
\begin{equation}
  \nabla_i U(q) + \lambda m_i q_i =0, \quad i=1,\dots,n. \label{eq:defCC}
\end{equation}
\end{definition}
A \emph{normalized CC} is a CC with $I(q)=1$.
Let us denote by NCC$_n$ the set of normalized central configurations for $n$-bodies.

The problem of infinite spin arises from the rotational symmetry of the $n$-body problem. If  $q \in \text{NCC}_n$, then   $R(\theta)q=(R(\theta)q_1,\dots,R(\theta)q_n) \in \text{NCC}_n$, where $R(\theta)$
is a rotation by angle $\theta$ of the plane. Therefore NCC$_n$ contains circles of similar CCs.

A well-known conjecture (see Smale \cite{Smale98}) is that, for a fixed choice of masses $m_i$, there are only finitely many CCs, up to symmetries of the space: orthogonal transformations, translations and scaling. In the planar case, it is a classical result that the conjecture holds for \(n=3\). More recent works (see \cite{AlbouyKaloshin12,HamtonMoeckel06}) have shown that it also holds for \(n=4\) and for generic masses when \(n=5\).

We call a planar normalized central configuration $q\in \text{NCC}_n$ \emph{isolated} if the equivalence class of $q$, denoted by $[q]$, is isolated in the quotient space $\text{NCC}_n/SO(2)$. We also call a configuration $q \in \text{NCC}_n$ \emph{non-degenerate} if $[q]$ is not a degenerate solution of system (\ref{eq:defCC}) in the space $S_n/SO(2)$.

Now, if $q(t)$ is a partial collision solution of \eqref{eq:second-order-formulation}, then for the colliding cluster $\mathcal{G}$ we can form curves $\hat{q}_\mathcal{G}(t) \in S_{\#\mathcal{G}}$ and
$\left[\hat{q}_\mathcal{G}(t)\right] \in S_{\#\mathcal{G}}/SO(2)$.
It follows from Chazy results \cite{Ch18}, generalized by Saari \cite{Saari71}, that $\left[\hat{q}_\mathcal{G}(t)\right]$ converges to a compact subset of  $\text{NCC}_{\#\mathcal{G}}/SO(2)$ as $t \to T$.  If $q_0$ is isolated  in $\text{NCC}_{\#\mathcal{G}}$ and  is in the limit set of $\hat{q}_\mathcal{G}(t)$, then   $\left[\hat{q}_\mathcal{G}(t)\right]$ converges
to $[q_0]$ in $S_{\#\mathcal{G}}/SO(2)$. The question of infinite spin is: whether $\hat{q}_\mathcal{G}(t)$ converges to a particular NCC or to a nontrivial subset of the circle $R(\theta)q_0$.  The main goal of this paper is to show that the latter cannot happen.  This generalizes the result of Moeckel and Montgomery \cite{Moeckel2023} for total collisions.

\begin{theorem}\label{thm:finite_spin}
    Suppose that $q(t)$ is a solution of $n$-body problem \eqref{eq:second-order-formulation} undergoing a partial collision as $t \to T$. Let $\mathcal{G}$ be a nontrivial cluster of colliding bodies. Suppose that
      its reduced and normalized configuration $\left[\hat{q}_\mathcal{G}(t) \right] \in S_{\#\mathcal{G}}/SO(2)$ converges to an isolated $\text{NCC}_{\#\mathcal{G}}$. Then  $\hat{q}_\mathcal{G}(t)$ converges to a particular CC.
\end{theorem}
The above result is claimed to be proven in \cite{E} even for the spatial case.  As explained in Appendix in  \cite{Moeckel2023} this proof (and also proof from \cite{SH} for the total collision) contains a gap related
to the `falling cat' phenomenon.

Just as in the case of total collision the problem is really to handle the degenerate case, because in the non-degenerate case the above theorem is proven by Chazy in \cite{Ch18} for total collisions while for partial collision the proof of the theorem for the non-degenerate case was outlined by Saari \cite{Saari84}. 
On the other hand, the degenerate but isolated CCs do exist for $n \geq 4$, see for example \cite{Palmore76} or literature cited in \cite{Moeckel2023}.

Our approach extends the one of \cite{Moeckel2023}. Let us briefly recall what is done there.  First, a good coordinate system is introduced in which the action of the rotation group is particularly simple and allows the separation of the angular variable $\theta$ from the others. The equations of motion (\ref{eq:second-order-formulation}) are then regularized
using McGehee coordinates and are of the following form
\begin{eqnarray}
  (r,v,s,b)'&=&f(r,v,s,b),  \label{eq:MM1} \\
  \theta'&=& g(s,b), \label{eq:MM-theta}
\end{eqnarray}
where the derivative is taken with respect to a new $\tau$ variable for which the collision occurs with $\tau\to \infty$, $r$ is the size of the configuration, $v$ is a scaled derivative of $r$, $s$ describe the shape of configuration and $b$ is its scaled derivative.
When the total collision happens then $(r,v,s,b)(\tau)$ for $\tau \to \infty$ approach a manifold of fixed points of the equation (\ref{eq:MM1}) contained in the set $\{(r=0,v,s,b=0), s\in \textrm{CC}\}$
where $\textrm{CC}$ is the set of scaled central configurations.  It is assumed in \cite{Moeckel2023} that $s(\tau) \to s_0$, which is an isolated shape in $\textrm{CC}$ and the question
of infinite spin is reduced showing that the length of the orbit $(r,v,s,b)(\tau)$ is finite.  This is achieved in \cite{Moeckel2023} by using a center manifold
near the limit point, in two steps: first, it is shown that the total collision orbit is approaching exponentially fast some orbit  $z_s(\tau)$ on the center manifold and
next using arguments based on the \L ojasiewicz inequality it is shown that orbit $z_s(\tau)$ has finite length, which implies finite length for the total collapse orbit.

In the derivation of system (\ref{eq:MM1},\ref{eq:MM-theta}) vanishing of angular momentum $\mu$ was crucial and the center of mass was fixed.

Now when considering a partial collision of a cluster of bodies $\mathcal{G}$, we can write similar equations for bodies in $\mathcal{G}$ and treat the influence of other bodies as a perturbation.
But we no longer have vanishing angular momentum of the bodies in $\mathcal{G}$ and the center of mass of cluster $\mathcal{G}$ is moving. We prove  that in the coordinate frame in which
the collision  point for the colliding cluster is fixed, the angular momentum (with respect to collision point) and its derivative decay exponentially fast with respect to the new time variable $\tau$
 and in the coordinates of Moeckel and Montgomery \cite{Moeckel2023} we obtain the following non-autonomous perturbation of equations (\ref{eq:MM1},\ref{eq:MM-theta})
 \begin{eqnarray}
  (r,v,s,b)'&=&f(r,v,s,b) + p_1(\tau,r,v,s,b,\mu,\mu'),  \label{eq:MM1-per} \\
  \theta'&=& g(s,b)+  p_2(\tau,r,v,s,b,\mu,\mu'), \label{eq:MM-theta-per}
\end{eqnarray}
where $p_{1,2}(\dots)$ decay exponentially fast with respect to $\tau$.

Moreover, just as in the total collision case the orbit $(r,v,s,b)(\tau)$ for $\tau \to \infty$ is approaching a manifold of fixed points of equation (\ref{eq:MM1}) contained in set $\{(r=0,v,s,b=0), s\in \textrm{CC}\}$.  So, in the total collision case, we assumed as in \cite{Moeckel2023} that $s(\tau) \to s_0$, which is an isolated shape in $\textrm{CC}$.
These developments are contained in  Section~\ref{sec:MM-cord} with tedious details relegated to Appendix~\ref{sec:dotomega}.

Next, we show that our partial collision orbit for bodies in $\mathcal{G}$ projected on variables $(r,v,s,b)$ is approaching exponentially fast some orbit  $z_1(\tau)$ of equation  (\ref{eq:MM1}). This is a content of Section~\ref{sec:shadow-isoseg}.  Then, in Section~\ref{sec:no-inf-spin}, we continue with the finite length argument for $z_1(\tau)$ as in \cite{Moeckel2023}. It turns out that also  contribution of term $p_2$ to the spin is finite
and we obtain our conclusion.

To realize the above argument in the paper we rely on asymptotic rates of various quantities for a cluster of bodies going to the collision. This is needed to establish
bounds on terms $p_1$ and $p_2$ in equations (\ref{eq:MM1-per},\ref{eq:MM-theta-per}). For this purpose, we adapted to the setting of partial collisions some results from \cite{Wintner41,PS68} about total collision.  This is a content of Appendix~\ref{sec:est-coll}.

In Appendix~\ref{sec:conv-cc} we give a proof of the convergence of normalized configurations to the set of CCs for the partial collision case. Finally, in Appendix~\ref{sec:TaubThm} we state two Tauberian theorems from Wintner's book \cite{Wintner41}, which are used in the derivation of the asymptotic bounds for a cluster of bodies going to the collision.


\section{Derivation of Moeckel and Montgomery's equations and McGehee regularization}
\label{sec:MM-cord}

Let us single out one cluster $\mathcal{G}$ with  $1<\# \mathcal{G} =k < n$ and focus on a \emph{collision orbit} $q_i$, where $i\in \mathcal{G}$. If we separate the right-hand side of its Newton equation \eqref{eq:second-order-formulation} in two terms, one with indexes in the cluster $\mathcal{G}$ and the other with indexes in other clusters, we get
\begin{align*}
    m_i \ddot{q}_i = \sum\limits_{i < j,\, j\in \mathcal{G} }\frac{m_im_j}{q_{ij}^3}(q_j - q_i) + \sum_{j\notin \mathcal{G}}\frac{m_im_j}{q_{ij}^3}(q_j - q_i).
\end{align*}
Since the bodies in the same cluster collide, we have a singularity in the first term.

The second-order differential equation \eqref{eq:second-order-formulation} for all bodies in cluster $\mathcal{G}$ is equivalent to the Lagrangian equations for the following Lagrangian function,
\begin{align}
    \mathcal{L}_\text{tot} = K(\dot{q}) + U(q) = \sum_{i=1}^n \frac{m_i\dot{q}_i^2}2  + U(q)\text{,} \nonumber
\end{align}
where the kinetic energy $K(\dot{q})=\frac{1}{2}I(\dot{q})$.

Studying a partial collision of bodies in the cluster $\mathcal{G}$, we will also investigate the part of the Lagrangian related to this cluster, that is,
\begin{align}\label{eq:class_lagrangian}
    \mathcal{L} :=\mathcal{L}_{\mathcal{G}} := K_{\mathcal{G}}(\dot{q}) + U_{\mathcal{G}}(q) + U_{\text{ext}}(q),
\end{align}
where $K_{\mathcal{G}}$ and $U_{\mathcal{G}}$ are, respectively, the kinetic energy and the potential function of the bodies restricted to bodies with indexes in $\mathcal{G}$
and $U_{\text{ext}}$ is the potential generated by the interaction of bodies inside and outside the cluster:
\begin{equation}\label{eq:energy_potential_inCluster}
K_{\mathcal{G}}(\dot{q}) = \sum_{i\in \mathcal{G}} \frac{m_i\dot{q}_i^2}2
\text{, \qquad }
U_{\mathcal{G}}(q) = \sum_{i<j;i,j\in \mathcal{G}}\frac{m_i m_j}{q_{ij}}
\text{, \qquad }
U_{\text{ext}}(q) = \sum_{i\in \mathcal{G}, j\notin \mathcal{G}}\frac{m_i m_j}{q_{ij}}.
\end{equation}

The equations of motion of orbits tending to the partial collision at $L_{\mathcal{G}}$ are the Lagrangian equations for $\mathcal{L}$.

A common approach to develop a Lagrangian of an $n$-body problem is to fix the center of mass $c$ at the origin and introduce the Jacobi coordinates to separate the terms describing the motion of $c$ and to reduce the number of variables. This method is also convenient in \cite{Moeckel2023}, in the case of total collision, the collision point would also be the origin. In our case of partial collision of the cluster $\mathcal{G}$, however,  $c_\mathcal{G}$ is moving and cannot be removed by choosing inertial coordinate frame. 

From now on, suppose that the system \eqref{eq:second-order-formulation} admits a partial collision of cluster $\mathcal{G}$ at $L_{\mathcal{G}}=0$ and, by reordering the indexes, let $\mathcal{G} = \{1, 2, \dots, k\}$ for $k<n$. If that would not lead to misunderstanding, we shall use the same notation $q = q(t)$ also for the position vector of the bodies in the cluster $\mathcal{G}$, that is $q(t)\in\mathbb{R}^{2k}$.

Let us define
\begin{eqnarray}
  I_\mathcal{G}(q)&=&\sum_{i \in \mathcal{G}} m_i q_i^2, \\
  M_{\mathcal{G}} &=& \sum\limits_{i\in \mathcal{G}}m_i.
\end{eqnarray}
Observe that $K_\mathcal{G}(\dot{q})=\frac{1}{2} I_\mathcal{G}(\dot{q})$ and
\begin{equation}
   I_\mathcal{G}(q)=I^0_\mathcal{G}(q) + M_\mathcal{G}c^2_{\mathcal{G}}. \label{eq:IIGcg2}
\end{equation}

Denote by $P: \mathbb{R}^{2k} \rightarrow \mathbb{R}^{2k}$ the linear transformation matrix of $q = Pz$ such that $z$ is the generalized Jacobian coordinates. In this coordinate system, we have that
\begin{eqnarray}
    z_k &=& \frac{1}{M_{\mathcal{G}}}\sum _{i=1}^k m_iq_i=c_\mathcal{G},   \label{eq:zk}\\
    I_\mathcal{G}(z)&=& z^T M z =  \fz^T\tilde{M}\fz + M_\mathcal{G}c_{\mathcal{G}}^2, \label{eq:IGz} \\
     U_\mathcal{G}(z) &=& U(\fz).  \label{eq:Ufz}
\end{eqnarray}
where $\mathfrak{z} = (z_1,\dots, z_{k-1})$, $M=P^T\text{diag}(m_1,m_1, \dots , m_k, m_k) P$ is a $2k\times 2k$  positive-definite symmetric matrix,
 and $\tilde{M}$  is a   positive-definite symmetric of size  $(2k-2)\times (2k-2)$.

The meaning of (\ref{eq:Ufz}) is as follows:  in new variables potential $U_\mathcal{G}$ depends on $\fz$ and does not depend on $z_k=c_\mathcal{G}$. It also defines function $U(\fz)$.

Moreover, in new coordinates the action of rotation group $SO(2)$, is still simple, namely for $R(\alpha)$ being rotation by $\alpha$ we have
\begin{equation}
R(\alpha)(z_1,\dots,z_k)=(R(\alpha)z_1,\dots,R(\alpha)z_k).
\end{equation}
And finally the collision point is given by $\fz=0,c_\mathcal{G}=0$.

Observe that
\begin{eqnarray}
  I^0_\mathcal{G}(q)&=&\fz^T\tilde{M}\fz=:I^0_\mathcal{G}(\fz), \label{eq:IIGcg1}
\end{eqnarray}
In the above expression we defined function $I^0_\mathcal{G}(\fz)$ by the same symbol $I^0_\mathcal{G}(q)$, but while their formulas differ their values are the same if
arguments $q$ and $\fz$ are related by coordinate transformation $P$.

Since $\dot{q} = P\dot{z}$, from (\ref{eq:IGz}) we have that the kinetic energy in \eqref{eq:class_lagrangian} becomes
\begin{align}
    K_{\mathcal{G}} = \frac{1}{2}\dot{\fz}^T\tilde{M}\dot{\fz} + \frac{M_\mathcal{G}}{2}\dot{c}_{\mathcal{G}}^2, \label{eq:cluster_kinetic_energy_MMcoord}
\end{align}

Therefore, the Lagrangian \eqref{eq:class_lagrangian} becomes
\begin{align}\label{eq:expanded_lagrangian}
    \mathcal{L} =  \frac{1}{2}\dot{\fz}^T\tilde{M}\dot{\fz} +U(\fz) + \frac{M_\mathcal{G}}{2}\dot{c}_{\mathcal{G}}^2  + U_{\text{ext}}.
\end{align} 
We perform a coordinate change as \cite{Moeckel2023} for the configuration $\fz$. First, define
\begin{equation}
r = \left\|\fz\right\| = \sqrt{\fz^T\tilde{M}\fz}= \sqrt{I^0_\mathcal{G}(\fz)} \nonumber 
\end{equation}
and the normalized configuration
\begin{equation*}
\hat{\fz} = \fz/r \in \mathbb{S}^{2k-3}.
\end{equation*}
As in \cite{Moeckel2023}, to introduce the angular coordinate, we consider the open set of $\mathbb{S}^{2k-3}$ where $\hat{\fz}_{k-1}\neq 0$.

For this it is convenient to view $z_j$'s and $\dot{z}_j$'s as complex numbers, so we have $\mathfrak{z} \in \mathbb{C}^{k-1}$.
Observe that in those coordinates the action of the $SO(2)$ group is just $z_j \to e^{i\alpha} z_j$. We introduce a Hermitian mass metric on $\mathbb{C}^{k-1}$
\begin{equation}
   \langle\langle u,v\rangle\rangle_{\mathbb{C}} = \bar{u}^T\tilde{M}v, \quad u,v\in\mathbb{C}^{k-1},  \label{eq:her-complex}
\end{equation}
where $\bar{u}$ is a conjugate of $u$. Observe that with the above definition $r=\sqrt{ \langle\langle \mathfrak{z},\mathfrak{z}\rangle\rangle_{\mathbb{C}}}$.

\begin{remark}
We note that \(\tilde{M}\) in (\ref{eq:her-complex}) is different from \(\tilde{M}\) in \eqref{eq:cluster_kinetic_energy_MMcoord}. Specifically, in the complex case, it is half the size of the matrix used in the real case.
However, to avoid unnecessary notation, we will allow ourselves this slight abuse of notation.
\end{remark}
Consider that $\hat{\fz}_{k-1}$ lies on the $OY$ axis. This represents a local section of the action of the rotation group.
By taking $\hat{\fz}_{k-1} = \kappa e^{i\theta}$, $\kappa>0$ and defining $s_i := \hat{\fz}_i\hat{\fz}^{-1}_{k-1} \in \mathbb{C}$, we have $\hat{\fz}_i = \kappa e^{i\theta}s_i$, that is,
\begin{align*}
    \hat{\fz} = \kappa e^{i\theta}(s,1)\text{,}
\end{align*}
where $s\in \mathbb{C}^{k-2}$. Since $\left\|\hat{\fz}\right\| =1$ and $\kappa = \left\|(s,1)\right\|^{-1}$, we have the following coordinate change
\begin{align}\label{eq:change_of_coordinates}
    \fz = re^{i\theta}\frac{(s,1)}{\left\|(s,1)\right\|}, \quad (r,\theta, s)\in \mathbb{R}^+\times \mathbb{S}^1\times \mathbb{C}^{k-2}.
\end{align}

\begin{remark}
    Note that the collision happens when $r=0$.
\end{remark}

\begin{lemma}
\label{lem:MM-kin-part}
Define $\rho = \dot{r}$ and $\omega = \dot{s}$.
Using the change of coordinates \eqref{eq:change_of_coordinates}, we have that
\begin{align}\label{eq:Kin-cluster}
    \frac{1}{2}\dot{\fz}^T\tilde{M}\dot{\fz}  = \frac{\rho^2}{2} + \frac{r^2\dot{\theta}^2}{2} +\frac{r^2\left\|\omega\right\|^2}{2\left\|(s,1)\right\|^2} -\frac{r^2 G(s,\omega)^2}{2\left\|(s,1)\right\|^4} + \frac{r^2\dot{\theta}\Omega(s,w)}{\left\|(s,1)\right\|^2},
\end{align}
    where $\left\|\omega\right\| = \left\|(\omega,0)\right\|$,  $G(s,\omega) = \operatorname{Re}\langle\langle (s,1), (\omega,0)\rangle\rangle_{\mathbb{C}}$ and $\Omega(s,\omega)= \operatorname{Im}\langle\langle (s,1), (\omega,0)\rangle\rangle_{\mathbb{C}}$.
\end{lemma}
\begin{proof}
    This equality is obtained from the same computations as in \cite{Moeckel2023}.
\end{proof}

Now, if we pass back to real coordinates $s,\omega \in \mathbb{R}^{2(k-2)}$, from \eqref{eq:Kin-cluster} we obtain the following expression for the kinetic energy term for $\mathfrak{z}$:
\begin{equation}
   \frac{1}{2}\dot{\fz}^T\tilde{M}\dot{\fz}  = \frac{\rho^2}{2} + \frac{r^2\dot{\theta}^2}{2} +\frac{r^2}{2}F(s,\omega) +\frac{r^2 \Omega(s,\omega)^2}{2\left\|(s,1)\right\|^4} + \frac{r^2\dot{\theta}\Omega(s,w)}{\left\|(s,1)\right\|^2},  \label{eq:KinPart-cluster}
\end{equation}
where $F(s,\omega)$ is the local coordinate representation of the square of the Fubini-Study norm (a K\"ahler metric on a complex projective space $\mathbb{C}P^n$ endowed with a Hermitian form) and takes the form
\begin{equation}\label{eq:Fubini-study_norm}
  F(s,\omega) = \omega^T\mathcal{A}(s)\omega =   \frac{\left\|\omega\right\|^2}{\left\|(s,1)\right\|^2} - \frac{G(s,\omega)^2 +\Omega(s,\omega)^2}{\left\|(s,1)\right\|^4},
\end{equation}
where $\mathcal{A}(s)$ is a positive definite $(2k-4)\times (2k-4)$ matrix. Besides, we remark that \\ $\Omega(s,\omega)~=~\mathcal{B}(s)^T\omega$, where $\mathcal{B}(s) \in \mathbb{R}^{2(k-2)}$ is linear in $s$ plus a constant term. 
\subsection{Equations of motion}
\label{secsub:eq-motion}

We now derive the equations of motion for bodies in cluster $\mathcal{G}$ via Lagrangian formulation.
First, we plug expression \eqref{eq:KinPart-cluster} into \eqref{eq:expanded_lagrangian}, then our Lagrangian function is
\begin{align}
    \mathcal{L} &= \frac{\rho^2}{2} + \frac{r^2\dot{\theta}^2}{2} +\frac{r^2}{2}F(s,\omega) +\frac{r^2 \Omega(s,\omega)^2}{2\left\|(s,1)\right\|^4} + \frac{r^2\dot{\theta}\Omega(s,w)}{\left\|(s,1)\right\|^2}+ \frac{1}{r}V(s) + U_{\text{ext}}+\frac{M_\mathcal{G}}{2}\dot{c}_{\mathcal{G}}^2,  \label{eq:LagMMcoord}
\end{align}
where
\begin{equation*}
  V(s) = \left\|(s,1)\right\|U(s,1).
\end{equation*}
The Euler-Lagrange equation for the motion of $\fz$ is
\begin{align*}
    \frac{d}{dt}\frac{\partial \mathcal{L}}{\partial \dot{\mathfrak{z}}} - \frac{\partial \mathcal{L}}{\partial \fz} = 0.
\end{align*}
However, we want to obtain those equations on  $(r,\rho,\theta,\dot{\theta},s,\omega)$ variables.
First, we define $\mu(t)$ as
\begin{align}
    \mu(t) := \partial_{\dot{\theta}}\mathcal{L} = r^2\dot{\theta} + \frac{r^2\Omega(s,\omega)}{\left\|(s,1)\right\|^2} \quad\quad \Rightarrow \quad\quad\dot{\theta} = \frac{\mu(t)}{r^2} - \frac{\Omega(s,\omega)}{\left\|(s,1)\right\|^2}.  \label{eq:mu-dot-theta}
\end{align}

We would like later use the asymptotic estimates for bodies going to collision to obtain asymptotic of $\mu$ and
$\dot{\mu}$. These asymptotic bounds are obtained in cartesian coordinates.  Therefore we need to have an expression for $\mu$ in cartesian coordinates.

\begin{remark}[Angular momentum in various coordinates]
Assume that we have a Lagrangian $L(\dot{w},w)$ and transformation of coordinates $w=T(q)$, then we have $\dot{w}=DT(q) \dot{q}$ and the Lagrangian in $q$-coordinates becomes
\begin{eqnarray*}
  \tilde{L}(\dot{q},q)=L(DT(q) \dot{q}, T(q)).
\end{eqnarray*}

Assume that we have a one-parameter group induced by an ODE
\begin{equation*}
   q'=\tilde{f}(q), \quad w'=f(w).
\end{equation*}
The relation between $\tilde{f}(q)$ and $f(w)$ is
\begin{eqnarray*}
   f(w)= w'= DT(q) q'=DT(q) \tilde{f}(q).
\end{eqnarray*}

We will show the following
\begin{lemma}
Under the above assumptions,
\begin{equation*}
  \frac{\partial \tilde{L}}{\partial \dot{q}}(\dot{q},q) q'(q) = \frac{\partial L}{\partial \dot{w}}(\dot{w},w) w'(w),
\end{equation*}
where $(\dot{w},w)=(DT(q)\dot{q},T(q))$.
\end{lemma}
\begin{proof}
We have
\begin{align*}
  \frac{\partial \tilde{L}}{\partial \dot{q}}(\dot{q},q) q'(q)&= \left(\frac{\partial L}{\partial \dot{w}}(DT(q)\dot{q},T(q)) \cdot \frac{\partial \dot{w}}{\partial \dot{q}}(q) \right) \cdot \tilde{f}(q)
   = \frac{\partial L}{\partial \dot{w}}(DT(q)\dot{q},T(q)) \cdot \left(DT(q)  \cdot \tilde{f}(q)\right)\\ &= \frac{\partial L}{\partial \dot{w}}(\dot{w},w) \cdot  f(w)= \frac{\partial L}{\partial \dot{w}}(\dot{w},w) w'(w).
\end{align*}
\end{proof}

 Let $R(\alpha)$ be the rotation by angle
$\alpha$  and let $J=R'(0)=\left[\begin{smallmatrix} 0 & -1 \\ 1 & 0 \end{smallmatrix} \right]$.

In Cartesian coordinates, the action of the rotation group is given by $q_i \mapsto R(\alpha) q_i$. Then,
\begin{equation*}
  q_i'=(-y_i,x_i)
\end{equation*}
and for the mechanical Lagrangian with the kinetic part given by $\sum_i \frac{m_i}{2} \dot{q}_i^2$, the \emph{angular momentum} of the bodies in $\mathcal{G}$ is given by
\begin{equation}
  \mu_\mathcal{G}=\frac{\partial L}{\partial \dot{q}}q'=\sum_{i\in \mathcal{G}} m_i (\dot{x}_i,\dot{y_i}) \cdot (-y_i,x_i)=\sum_{i\in\mathcal{G}} m_i q_i \times \dot{q}_i \label{eq:ang-mom-cart}.
\end{equation}
It is important to remark that in this context, the cross product `$\times$' gives a scalar. 
Since the vectors are in the plane, it takes the same value as the third coordinate of standard vector product `$\times$' in $\mathbb{R}^3$.

On the other side, the action of $SO(2)$ in new coordinates is  $(r,\theta,s,c_{\mathcal{G}}) \mapsto (r,\theta+\alpha,s,R(\alpha)c_{\mathcal{G}}) $, hence $(r,\theta,s,c_{\mathcal{G}})'=(0,1,0,Jc_{\mathcal{G}})$ and so we obtain that
\begin{equation*}
  \mu_\mathcal{G}=\frac{\partial L}{\partial \dot{\theta}} + \frac{\partial L}{\partial \dot{c}_{\mathcal{G}}} c_{\mathcal{G}}'= \frac{\partial L}{\partial \dot{\theta}} + M_\mathcal{G} c_{\mathcal{G}} \times \dot{c}_{\mathcal{G}}=\mu + M_\mathcal{G} c_{\mathcal{G}} \times \dot{c}_{\mathcal{G}}.
\end{equation*}

Therefore,
\begin{equation}
  \mu= \sum_{i\in\mathcal{G}} m_i q_i \times \dot{q}_i -  M_{\mathcal{G}} c_{\mathcal{G}}  \times \dot{c}_{\mathcal{G}}.   \label{eq:mu-corr-cart}
\end{equation}

It is easy to see that $\mu$ is equal to intrinsic angular momentum of cluster $\mathcal{G}$ as defined in  \cite{E} and \cite{Sp} (we use our notation) by
\begin{equation*}
  \mu^0_{\mathcal{G}}=\sum_{i \in \mathcal{G}} m_i (q_i - c_\mathcal{G})\times \left(\frac{d}{dt}(q_i - c_\mathcal{G})\right).
\end{equation*}

\end{remark}

\begin{remark}
For a system with the center of mass fixed at $0$, $\mu$ represents the angular momentum. It is known (see for example \cite{S13,Wintner41})  a total collapse collision only happens for $\mu = 0$ .
However, the same is not true for partial collisions.
\end{remark}

As an immediate consequence of definition of $\mu$ in \eqref{eq:mu-dot-theta}, the Euler-Lagrange equation for  $\theta$ and $\dot{\theta}$ is given by
\begin{equation*}
    \frac{d}{dt}\frac{\partial \mathcal{L}}{\partial \dot{\theta}}=  \frac{\partial \mathcal{L}}{\partial \theta},
\end{equation*}
hence
\begin{equation*}
    \dot{\mu} = \frac{\partial U_{\text{ext}}}{\partial \theta}.\nonumber 
\end{equation*}
Unlike the case of total collapse, the $\theta$ angle is not a cyclic variable for $\mathcal{L}$.

For variable $\rho$, the Euler-Lagrange equation takes the following form:
\begin{align}
    \frac{d}{dt}\frac{\partial \mathcal{L}}{\partial \rho} &=\frac{\partial \mathcal{L}}{\partial r}  \nonumber\\
    \dot{\rho} &= r\left(\dot{\theta}^2 + F(s,\omega) +\frac{\Omega(s,\omega)^2}{\left\|(s,1)\right\|^4} + \frac{2\dot{\theta}\Omega(s,w)}{\left\|(s,1)\right\|^2}\right) - \frac{V(s)}{r^2}  + \frac{\partial U_{\text{ext}}}{\partial r},
\end{align}
hence after using  (\ref{eq:mu-dot-theta}) we obtain
\begin{equation*}
    \dot{\rho} = rF(s,\omega) - \frac{V(s)}{r^2} + \frac{\mu^2}{r^3} + \frac{\partial U_{\text{ext}}}{\partial r}.  \nonumber 
\end{equation*}

For variable $\omega$ (for derivation see Appendix~\ref{sec:dotomega}):
\begin{eqnarray*}
    \frac{d}{dt}\frac{\partial \mathcal{L}}{\partial \omega} &=&\frac{\partial \mathcal{L}}{\partial s} \nonumber\\
    \dot{\omega} &=& \frac{1}{r^3}\mathcal{A}^{-1}(s)\nabla V(s) -\frac{2\rho \omega}{r} +\frac{1}{2}\mathcal{A}^{-1}(s)\nabla F(s,\omega) - \mathcal{A}^{-1}(s)(\nabla_s \mathcal{A}(s)\omega)\omega \nonumber\\
     & & + \frac{1}{r}\frac{\partial U_{\text{ext}}}{\partial s} -  \frac{\dot{\mu}}{r^2}\mathcal{A}(s)^{-1}\frac{\mathcal{B}(s)}{\|(s,1)\|^2}.\label{eq:lagrangeEquationW}
\end{eqnarray*}

Therefore, we have the following lemma 
\begin{lemma}
\label{lem:lagEqgoodCoord}
    In the variables $r,\rho, \theta,\mu, s$ and $\omega$, we have the following equations of motion of cluster~$\mathcal{G}$:
    \begin{align*}
    &\dot{r} = \rho\\
        &\dot{\rho} = rF(s,\omega) - \frac{V(s)}{r^2} + \frac{\mu^2}{r^3} + \frac{\partial U_{\text{\rm ext}}}{\partial r} .\\
        &\dot{s} = \omega \\
        &\dot{\omega} = \frac{1}{r^3}\mathcal{A}^{-1}(s)\nabla V(s) -\frac{2\rho \omega}{r} +\frac{1}{2}\mathcal{A}^{-1}(s)\nabla F(s,\omega) - \mathcal{A}^{-1}(s)(\nabla_s \mathcal{A}(s)\omega)\omega \nonumber\\
     &  + \frac{1}{r}\frac{\partial U_{\text{\rm ext}}}{\partial s} -  \frac{\dot{\mu}}{r^2}\mathcal{A}(s)^{-1}\frac{\mathcal{B}(s)}{\|(s,1)\|^2}\\
        &\dot{\theta} = \frac{\mu(t)}{r^2} - \frac{\Omega(s,\omega)}{\left\|(s,1)\right\|^2} \\
        & \dot{\mu} = \frac{\partial U_{\text{\rm ext}}}{\partial \theta}.
    \end{align*}
\end{lemma}

In the sequel, we will concentrate on the subsystem consisting of the first four variables ($r,\rho,s,\omega$). Such a system can be regarded as a non-autonomous system, with terms   $\frac{\mu^2}{r^3}+\frac{\partial U_{\text{ext}}}{\partial r}$, $\frac{1}{r}\frac{\partial U_{\text{ext}}}{\partial s} +\frac{\dot{\mu}}{r^2}\mathcal{A}(s)^{-1}\frac{\mathcal{B}(s)}{\|(s,1)\|^2}$  treated as non-autonomous perturbations.

For the cluster $\mathcal{G}$, we have energy $H_\mathcal{G}=K_\mathcal{G} - U_\mathcal{G}$. In the new coordinates, since $\dot{\theta} = \mu/r^2 - \Omega(s,\omega)/\left\|(s,1)\right\|^2$, this energy becomes
\begin{eqnarray}\label{eq:H-MMcoord}
  H_\mathcal{G} = \frac{\rho^2}{2} + \frac{\mu^2}{2r^2} +\frac{r^2}{2}F(s,\omega)  -\frac{1}{r}V(s) + \frac{M_\mathcal{G}}{2}\dot{c}_{\mathcal{G}}^2.
\end{eqnarray}

\subsection{McGehee-like regularization}
\label{subsec:mcgehee-reg}

The system of equations from Lemma~\ref{lem:lagEqgoodCoord}  exhibits a singularity at $r=0$. Then, to study these equations close to the collision manifold, it is necessary to remove this singularity.
In this order, McGehee's blow-up technique is used \cite{McGehee74, McGehee78}:
\begin{align}
    v = \sqrt{r}\rho, \quad \zw = r^{\frac{3}{2}}\omega  \label{eq:McGhee}
\end{align}
and  a new time variable $\tau$ such that $()' = r^{\frac{3}{2}} \dot{()}$.

\begin{remark}
In contrast to \cite{Moeckel2023}, we use the variable $\zw$ instead of $w$, which looks similar to $\omega$ in print and may easily cause confusion.
\end{remark}

In these new variables, the differential equations for the variables $(r,v, s, \zw)$ are
\begin{align}
    r'  &= vr \label{eq:derivative of r-McGehee} \\
    v' &= \frac{v^2}{2} + F(s, \zw) -V(s) + \frac{\mu^2}{r} + r^2\frac{\partial U_{\text{ext}}}{\partial r} \label{eq:v'}\\
    s' & = \zw \label{eq:s'}\\
    \zw' & = -\frac{v\zw}{2} + \tilde{\nabla}V(s) + \frac{1}{2}\tilde{\nabla}F(s,\zw) - \mathcal{A}^{-1}(s)\left(D\mathcal{A}(s)\zw\right)\zw + r^2\frac{\partial U_{\text{ext}}}{\partial s}
    -  r\dot{\mu}\frac{\mathcal{A}(s)^{-1} \mathcal{B}(s)}{\|(s,1)\|^2} \label{eq:derivativeOfOmega}
\end{align}
where $\tilde{\nabla} = \mathcal{A}^{-1}\nabla$.

\begin{remark}
    In \cite{Moeckel2023}, the authors defined the covariant derivative (related to the Fubini-Study metrics) $\tilde{D}_{\tau}\zw$  as
    \begin{equation*}
        \tilde{D}_{\tau}\zw := \zw' -  \frac{1}{2}\tilde{\nabla}F(s,\zw) + \mathcal{A}^{-1}(s)\left(D\mathcal{A}(s)\zw\right)\zw,
    \end{equation*}
    then Eq.\ \eqref{eq:derivativeOfOmega} can be rewritten as
    \begin{align*}
        \tilde{D}_{\tau}\zw = -\frac{v\zw}{2} +\tilde{\nabla}V(s) +  r^2\frac{\partial U_{\text{ext}}}{\partial s}
    -  r\dot{\mu}\frac{\mathcal{A}(s)^{-1} \mathcal{B}(s)}{\|(s,1)\|^2}.
    \end{align*}
\end{remark}

Note that after this change of variables, the differential equation for $\theta$ is
\begin{equation}
    \theta' = \frac{\mu}{\sqrt{r}} -\frac{\Omega(s,\zw)}{\left\|(s,1)\right\|^2}\label{eq:theta'}.
\end{equation} 
\subsection{Estimation for non-autonomous terms}
\label{subsec:estm-non-auto}

The goal of this subsection is to estimate the asymptotic behavior of the non-autonomous terms present in the system of equations  (\ref{eq:derivative of r-McGehee}--\ref{eq:derivativeOfOmega}) obtained
after McGehee regularization. These terms are
\begin{eqnarray}\label{eq:nonautonomous-perturbation-terms}
\frac{\mu^2}{r} + r^2\frac{\partial U_{\text{ext}}}{\partial r}
\quad \text{ and } \quad
r^2\frac{\partial U_{\text{ext}}}{\partial s}  -  r\dot{\mu}\frac{\mathcal{A}(s)^{-1} \mathcal{B}(s)}{\|(s,1)\|^2}.
\end{eqnarray}


For this subsection, we fix a partial collision orbit with trajectory $q(t)$ which starts at $t=t_0$ and goes to partial collision for $t \nearrow T$.
We change the time as in McGehee's transformation, then the relation between the new time variable $\tau$ and $t$ is
\begin{equation}
  d \tau = \frac{1}{r(t)^{3/2}}dt.  \label{eq:dtau/dt}
\end{equation}
We may set $\tau(t_0)=0$ and we have that $\tau(t)\to \infty$ as $t\nearrow T$. We shall prove that the terms \eqref{eq:nonautonomous-perturbation-terms} decay exponentially fast on a partial collision orbit as $\tau\to \infty$.
Our statements are based on general estimations of physical quantities near a partial collision singularity that are explored in Appendix \ref{sec:est-coll}, Theorems
 \ref{thm:r-estm-coll-time} and \ref{thm:intr-angmom-collison}.

We begin with some consequences of Theorem~\ref{thm:r-estm-coll-time}.
\begin{lemma}
\label{lem:asympt-bnds}
\begin{eqnarray}
  r(t) &\sim& A^{1/2}(T-t)^{2/3},   \label{eq:r(t)asympt}\\
  \rho(t)=\dot{r}(t) &\sim& -\frac{2}{3}A^{1/2}(T-t)^{-1/3}, \label{eq:dot-r(t)asympt} \\
  \ddot{r}(t)=\dot{\rho}(t) &\sim& -\frac{2}{9}A^{1/2}(T-t)^{-4/3},  \nonumber
\end{eqnarray}
where by $a(t) \sim b(t)$ we denote that\quad $\lim\limits_{t\to T}\frac{a(t)}{b(t)} \in (0,+\infty)$.
\end{lemma}
\begin{proof}
Since $r=\sqrt{I^0_\mathcal{G}}$  from (\ref{eq:r-thmJ-estm},\ref{eq:r-thmdotJ-estm},\ref{eq:r-ddotJ-bnds}) in Theorem~\ref{thm:r-estm-coll-time}  we obtain our assertion   after some easy computations.
\end{proof}

The following lemma, which describes the asymptotic behavior of $r$ with respect to $\tau$.
\begin{lemma}
\label{lem: rate convergence of r}
There exist positive constants $K_1$, $K_2$, $E_1$, $E_2$ such that the following inequalities hold for $\tau >0$:
\begin{equation}\label{eq:r(tau)exp}
  K_2e^{-E_2 \tau} \leq |r(\tau)| \leq K_1 e^{- E_1\tau},
\end{equation}
Moreover, for every $\varepsilon>0$ the constants may be chosen such that $E_2-E_1<\varepsilon$ and \eqref{eq:r(tau)exp} holds for $\tau > \tau_0\geq 0$.
\end{lemma}
\begin{proof}
From Lemma~\ref{lem:asympt-bnds}
\begin{equation}
\lim_{t \nearrow T} \frac{r(t)}{|T-t|^{\frac{2}{3}}} = \sqrt{A} \text{, \quad  for some constant $A>0$.} \label{eq:r-beha}\nonumber
\end{equation}
From the above estimate, it follows that there exist $B_1>A^{1/2}>B_2$ and $B_i\approx A^{1/2}$, such that  \
\begin{equation}
  B_2 (T-t)^{2/3} < r(t) <   B_1(T-t)^{2/3}.  \label{eq:Saari_r(t)}
\end{equation}
Therefore, we obtain the following relation from (\ref{eq:dtau/dt},\ref{eq:Saari_r(t)})
\begin{equation*}
 \frac{1}{B_1^{3/2} (T-t)} < \frac{d\tau}{dt} < \frac{1}{B_2^{3/2} (T-t)},
\end{equation*}
hence we obtain for $t \in (t_0,T)$
\begin{equation}
   B_1^{-3/2} \ln \frac{T-t_0}{T-t} <  \tau(t) - \tau(t_0)=\tau(t) < B_2^{-3/2} \ln \frac{T-t_0}{T-t}.  \label{eq:tau(t)}
\end{equation}
Observe that $\tau(t)\in [0,\infty)$  for $t \in [t_0,T)$.

Now we are ready to express $r$ in terms of $\tau$. From \eqref{eq:tau(t)}, it follows that
\begin{equation*}
  \exp\left(-\tau B_1^{\frac{3}{2}}\right)(T-t_0)  < T-t <  \exp\left(-\tau B_2^{\frac{3}{2}}\right)(T-t_0).
\end{equation*}
Hence from \eqref{eq:Saari_r(t)}, we obtain
\begin{equation*}
B_2 (T-t_0)^{2/3} \exp\left(-\frac{2\tau B_1^{\frac{3}{2}}}{3}\right)  <  r(\tau) < B_1 (T-t_0)^{2/3} \exp\left(-\frac{2\tau B_2^{\frac{3}{2}}}{3}\right).
\end{equation*}
\end{proof}

Since $U_{\text{ext}}$ are smooth functions, their partial derivatives are bounded. By Lemma~\ref{lem: rate convergence of r}, we have the following result.
\begin{lemma}
\label{lem:non-auto-terms}
    On a partial collision orbit, $r^2\frac{\partial U_{\text{\rm ext}}}{\partial r} = O(e^{-2E_1 \tau})$  and $r^2\frac{\partial U_{\text{\rm ext}}}{\partial s} = O(e^{-2E_1 \tau})$.
\end{lemma}


The last step is to estimate the asymptotic behavior of $\mu$, $\dot\mu$ and the terms containing them.
\begin{lemma}
\label{lem:mu2/r}
There exists a constant $E>0$ such that
    \begin{align*}
  \mu=O(r^{7/2}), \qquad \dot{\mu}=O(r^2), \quad\quad
       \frac{\mu(\tau)^2}{r(\tau)} = O(e^{-E \tau}), \qquad  r\dot{\mu}\frac{\mathcal{A}(s)^{-1} \mathcal{B}(s)}{\|(s,1)\|^2}=O(e^{-E \tau}).
    \end{align*}
\end{lemma}

\begin{proof}
From Theorems~\ref{thm:cmass-beha} and~\ref{thm:intr-angmom-collison}  we obtain $\mu(t)=O\left(|T-t|^{7/3}\right)$ and $\dot{\mu}(t)=O\left(|T-t|^{4/3}\right)$.
From the this and (\ref{eq:r(t)asympt}) in Lemma~\ref{lem:asympt-bnds}  we immediately conclude that $\mu = O(r^{7/2})$ and $\dot{\mu}=O(r^2)$. Therefore $ \mu^2/r= O(r^6)$ and
Lemma \ref{lem: rate convergence of r} ensures that
\begin{align*}
   \frac{\mu^2}{r}  = O(e^{-6E_1 \tau}).
\end{align*}

Finally, we also have
\begin{eqnarray*}
  r\dot{\mu}\frac{\mathcal{A}(s)^{-1} \mathcal{B}(s)}{\|(s,1)\|^2} = O(r^3)=O(e^{-3E_1 \tau}).
\end{eqnarray*}
\end{proof}

As a conclusion, the non-autonomous perturbation terms \eqref{eq:nonautonomous-perturbation-terms} converge exponentially fast to zero as $\tau \to \infty$, that is,
\begin{proposition}\label{prop:deltaEstimation}
    On a partial collision orbit, for $\tau$ large enough, there exist  positive constants $C$ and $E$ such that
    \begin{align*}
        \max\left\{\left\|\frac{\mu^2}{r} + r^2\frac{\partial U_{\text{\rm ext}}}{\partial r}\right\|, \left\|\quad r^2\frac{\partial U_{\text{\rm ext}}}{\partial s}  -  r\dot{\mu}\frac{\mathcal{A}(s)^{-1} \mathcal{B}(s)}{\|(s,1)\|^2}\right\|\right\} \leq Ce^{-E\tau}.
    \end{align*}
\end{proposition}
This proposition \ref{prop:deltaEstimation} will be essential in proving Thm.~\ref{thm:finite_spin}. 

\subsection{Convergence in the system after McGehee regularization}\label{sec:convergence}

The goal of this section is to  prove that our partial collision orbit converges to the set of CCs on the manifold of fixed points of system (\ref{eq:derivative of r-McGehee}--\ref{eq:derivativeOfOmega}) without non-autonomous terms, i.e., of the system
\begin{eqnarray}
     r'  &=& vr \label{eq:r'-auto} \\
    v' &=& \frac{v^2}{2} + F(s, \zw) -V(s)  \label{eq:v'-auto}\nonumber\\
    s' &=& \zw \label{eq:s'-auto} \nonumber\\
    \zw' &=& -\frac{v\zw}{2} + \tilde{\nabla}V(s) + \frac{1}{2}\tilde{\nabla}F(s,\zw) - \mathcal{A}^{-1}(s)\left(D\mathcal{A}(s)\zw\right)\zw.  \label{eq:zw'-auto}
\end{eqnarray}
In \cite{Moeckel2023}, this question for total collision was not discussed, but the reader can find it in another paper by Moeckel \cite{MoeckelTCSC}, also \cite{McGehee78,S13}.
\com{\\PZ: we should check what is in McGehee's papers about this.\\
AG: McGehee mentions it as a classical result of Sundman, particularly in 3bp\\}

\begin{lemma}
\label{lem:q/rtoncc}
$\hat{q}_\mathcal{G}(t)$ tends for $t \to T$ to the set of normalized central configurations.
\end{lemma}
\begin{proof}
We know from Theorem~\ref{thm:conver-to-cc}  that $\frac{q_\mathcal{G}-L_\mathcal{G}}{|T-t|^{2/3}}$  approaches the set  of central configurations. We have
\begin{eqnarray*}
  \hat{q}_\mathcal{G}(t)=\frac{q_\mathcal{G}(t)-c_\mathcal{G}(t)}{r(t)}= \frac{(q_\mathcal{G}(t)- L_\mathcal{G}) + (L_\mathcal{G}-c_\mathcal{G}(t))}{|T-t|^{2/3}} \frac{|T-t|^{2/3}}{r(t)}
\end{eqnarray*}
From Lemma~\ref{lem:asympt-bnds} we have that $\frac{|T-t|^{2/3}}{r(t)} \to A^{-1/2}$ and from Theorem~\ref{thm:cmass-beha} we know $L_\mathcal{G}-c_\mathcal{G}(t))=O(|T-t|)$,
therefore $\hat{q}_\mathcal{G}(t)$ also approaches the set of central configurations, which are normalized.
\end{proof}

\begin{lemma}
\label{lem:conv-to-man}
Assume that the limit set of $\hat{q}_\mathcal{G}(t)$ is contained in the domain, where coordinates from Section~\ref{sec:MM-cord} are defined.
Then the following holds   for $\tau \to \infty$
\begin{eqnarray*}
  r \to 0, \quad r' \to 0, \\
  v \to -\frac{2}{3} A^{3/4}, \quad v' \to 0, \\
  b \to 0, \quad b' \to 0.
\end{eqnarray*}
Moreover, $s$ is converging to the set $\{s:  \   \nabla V(s)=0\}$.
\end{lemma}
\noindent
\begin{proof}
From Lemma~\ref{lem:q/rtoncc} we know $\hat{q}_\mathcal{G}(t)$  approaches the set of normalized central configurations, hence it is away from singularities of $V$, see \cite{Sh70,MZ19}.
In our coordinates, such a set of normalized CC is given by condition $\nabla V(s)=0$.

From our assumption, it also follows that $\mathcal{A}(s)$ and its derivatives are bounded on our collision trajectory.

Obviously, $r \to 0$ because it is a collision orbit.
From Lemma~\ref{lem:asympt-bnds} and definition of $v$ in \eqref{eq:McGhee}, we obtain for $t \nearrow T$
\begin{eqnarray*}
  v(t)=\sqrt{r(t)} \rho(t) = \left(\sqrt{\frac{r(t)}{(T-t)^{2/3}}} (T-t)^{1/3}\right) \left(\frac{\rho(t)}{(T-t)^{-1/3}} (T-t)^{-1/3}\right)\\
    = \sqrt{\frac{r(t)}{(T-t)^{2/3}}} \frac{\rho(t)}{(T-t)^{-1/3}} \to A^{1/4} \left(-\frac{2}{3}A^{1/2}\right)=-\frac{2}{3} A^{3/4}.
\end{eqnarray*}
We see that $v$ approaches a fixed value. From this and (\ref{eq:derivative of r-McGehee}), it follows that $r' \to 0$.

Now we look at $v'$.
We have
\begin{equation*}
  v'(\tau)=\frac{dv}{dt} \cdot \frac{dt}{d\tau}=\dot{v} r^{3/2}.  \label{eq:v'dotv}
\end{equation*}
Observe that
\begin{eqnarray*}
  \dot{v}=\frac{d}{dt}(\sqrt{r}\rho)=\frac{\dot{r}}{2\sqrt{r}} \rho + \sqrt{r} \dot{\rho}=\frac{\rho^2}{2 \sqrt{r}} + \sqrt{r}\dot{\rho}.
\end{eqnarray*}

Now using Lemma~\ref{lem:asympt-bnds} we obtain
\begin{eqnarray*}
 v'=\dot{v}r^{3/2} = \frac{1}{2} \rho^2 r + r^2 \dot{\rho} \sim \frac{1}{2}\left(-\frac{2}{3}A^{1/2}(T-t)^{-1/3} \right)^2 \left(A^{1/2}(T-t)^{2/3}\right) \\
    + \left(A^{1/2}(T-t)^{2/3}\right)^2 \left( -\frac{2}{9}A^{1/2}(T-t)^{-4/3} \right)
    = \frac{2 A^{3/2}}{9} -  \frac{2 A^{3/2}}{9}=0.
\end{eqnarray*}
The mathematical interpretation of these computations is as follows: both terms $ \frac{1}{2} \rho^2 r$  and $r^2 \dot{\rho}$ have finite limits and these limits cancel out.
Therefore, we proved that $v' \to 0$.

By Theorem~\ref{thm:estm-coll-time},  energy $H_\mathcal{G}$, given by (\ref{eq:H-MMcoord}), is bounded.
 In McGehee coordinates, from (\ref{eq:H-MMcoord}) we obtain
\begin{equation*}
 \frac{1}{2}v^2 + \frac{1}{2}F(s,b) - V(s) + \frac{\mu^2}{2r} =r (H_\mathcal{G}(\tau) -M_\mathcal{G}\dot{c}_\mathcal{G}^2/2) .  \label{eq:ener-new-var}
\end{equation*}
Hence from (\ref{eq:v'}), Lemma~\ref{lem:mu2/r} and since $\frac{\partial U_{ext}}{\partial r}=O(1)$ we obtain
\begin{equation*}
  v'=\frac{1}{2}F(s,b) + \frac{\mu^2}{2r} + r^2 \frac{\partial U_{ext}}{\partial r}+ r(H_\mathcal{G}(\tau) -M_\mathcal{G}\dot{c}_\mathcal{G}^2/2)=
     \frac{1}{2}F(s,b) + O(r^2)+ r(H_\mathcal{G}(\tau) -M_\mathcal{G}\dot{c}_\mathcal{G}^2/2).
\end{equation*}
Since $v' \to 0$, $r \to 0$ and, by Theorems~\ref{thm:cmass-beha} and \ref{thm:estm-coll-time}, $H_\mathcal{G}(\tau)$ and $M_\mathcal{G}\dot{c}_\mathcal{G}^2/2$ are bounded, we have that $F(s,b) \to 0$, hence it is easy to see that from (\ref{eq:Fubini-study_norm}) it follows that $b \to 0$.

To prove that $b' \to 0$, we apply Theorem~\ref{thm:taubwin}. For this we have to show  that $b''(\tau)$ is bounded.

From (\ref{eq:derivativeOfOmega}) it follows that $b'$ is bounded because, by our assumptions, we are away from singularities of $V(s)$ and $\mathcal{A}(s)$, so these functions are bounded and non-autonomous terms converge to $0$ for $\tau \to \infty$ (see Lemmas~\ref{lem:non-auto-terms} and~\ref{lem:mu2/r}).
From (\ref{eq:derivativeOfOmega}) we have
\begin{eqnarray*}
  b''=  \frac{1}{2}\left(\tilde{\nabla}F(s,b)\right)' +\left( \mathcal{A}^{-1}(s) (D\mathcal{A}(s)(b))b\right)' + \left(\mathcal{A}(s)^{-1}\nabla V(s)\right)' - \frac{1}{2}\left(vb\right)' \\
   + \left(r^2\frac{\partial U_{\text{ext}}}{\partial s} \right)'
    -  \left(r\dot{\mu}\frac{\mathcal{A}(s)^{-1} \mathcal{B}(s)}{\|(s,1)\|^2}\right)'
\end{eqnarray*}
and
\begin{eqnarray*}
  \left(\tilde{\nabla}F(s,b)\right)'= \|\nabla_s \left((\mathcal{A}^{-1}(s) \nabla_s \mathcal{A}(s) \right) )\| \cdot |s'| \cdot |b|^2 + O(|b|\cdot |b'|) \to 0, \\
  \left( \mathcal{A}^{-1}(s) (D\mathcal{A}(s)(b))b\right)'=O(|b|) \to 0, \\
  \left(\mathcal{A}(s)^{-1}\nabla V(s)\right)'= O(s')=O(|b|) \to 0, \\
  \left(vb\right)'=v'b + vb' =O(1).
\end{eqnarray*}

Now we take a look at $\left(r^2\frac{\partial U_{\text{ext}}}{\partial s} \right)'$.  We  have
\begin{eqnarray*}
  \left(r^2\frac{\partial U_{\text{ext}}}{\partial s} \right)'=2rr'\frac{\partial U_{\text{ext}}}{\partial s} + r^2 \sum_{j \notin \mathcal{G}}\frac{\partial^2 U_{\text{ext}}}{\partial s \partial q_j} q_j' + r^2 \frac{\partial^2 U_{\text{ext}}}{\partial s^2}s'
\end{eqnarray*}
Terms $2r' r\frac{\partial U_{\text{ext}}}{\partial s}$ and $ r^2 \frac{\partial^2 U_{\text{ext}}}{\partial s^2}s'=r^2 \frac{\partial^2 U_{\text{ext}}}{\partial s^2}b$ are converging to $0$.

For the remaining term, since $q_j'=r^{3/2} \dot{q}_j$ we have
\begin{eqnarray*}
r^2 \sum_{j \notin \mathcal{G}}\frac{\partial^2 U_{\text{ext}}}{\partial s \partial q_j} q_j' =  r^{7/2} \sum_{j \notin \mathcal{G}}\frac{\partial^2 U_{\text{ext}}}{\partial s \partial q_j} \dot{q}_j.
\end{eqnarray*}
Observe that $\frac{\partial^2 U_{\text{ext}}}{\partial s \partial q_j}$ is bounded.

From Theorems~\ref{thm:r-estm-coll-time} and~\ref{thm:estm-coll-time}   applied to cluster $\mathcal{G}$ to get a bound on $r$ and to other clusters to get bound for $\dot{q}_j$ for $j \notin \mathcal{G}$ we obtain
\begin{eqnarray*}
   r^{7/2} \sum_{j \notin \mathcal{G}}\frac{\partial^2 U_{\text{ext}}}{\partial s \partial q_j} \dot{q}_j = O((T-t)^{\frac{2}{3} \cdot \frac{7}{2}}) \cdot O((T-t)^{-1/3})=
   O((T-t)^{2}).
\end{eqnarray*}

Let us denote $\mathcal{C}(s)=\frac{\mathcal{A}(s)^{-1} \mathcal{B}(s)}{\|(s,1)\|^2}$. Then $\mathcal{C}(s)$ is bounded and has bounded derivatives. We have
\begin{eqnarray*}
 \left(r\dot{\mu} \mathcal{C}(s)\right)'= r' \dot{\mu} \mathcal{C}(s) + r \dot{\mu}' \mathcal{C}(s) + r \dot{\mu} \mathcal{C}'(s)b.
\end{eqnarray*}

We have (we use Lemma~\ref{lem:mu2/r} to estimate $\dot{\mu}$) for $\tau \to \infty$
\begin{eqnarray*}
  r' \dot{\mu} \mathcal{C}(s) =  r' O(r) \to 0,  \\
   r \dot{\mu} \mathcal{C}'(s)b = O(r^2)b \to 0,
\end{eqnarray*}

Now we look at $r \dot{\mu}' \mathcal{C}(s)$. From Lemma~\ref{lem:lagEqgoodCoord} we have that $\dot{\mu}= \frac{\partial U_{\text{ext}}}{\partial \theta}$.  We express $\frac{\partial U_{\text{ext}}}{\partial \theta}$ in original cartesian coordinates.  For this observe that the partial derivative with respect to $\theta$ is the derivation with respect to
$\theta$  of $U_{\text{ext}}(R(\theta)q_1,R(\theta)q_2,\dots,R(\theta)q_k, q_{k+1},\dots,q_n)$  for $\theta=0$.  Therefore we obtain
\begin{equation*}
  \dot{\mu}=\sum_{i\leq k}\frac{\partial U_{\text{ext}}}{\partial q_i}Jq_i
\end{equation*}
where $J=R'(0)=\left[\begin{smallmatrix} 0 & -1 \\ 1 & 0 \end{smallmatrix} \right]$.

Hence
\begin{eqnarray*}
   \dot{\mu}'=\sum_{i \leq k} \sum_{j=1}^n  \frac{\partial^2 U_{\text{ext}}}{\partial q_i \partial q_j}Jq_i q_j' +  \sum_{i\leq k}\frac{\partial U_{\text{ext}}}{\partial q_i}Jq_i'= r^{3/2}\sum_{i \leq k} \sum_{j=1}^n  \frac{\partial^2 U_{\text{ext}}}{\partial q_i \partial q_j}Jq_i \dot{q}_j +  r^{3/2}\sum_{i\leq k}\frac{\partial U_{\text{ext}}}{\partial q_i}J\dot{q}_i
\end{eqnarray*}

The partial derivatives of $U_{\text{ext}}$ are bounded in the vicinity of collision.
From Theorems~\ref{thm:estm-coll-time} and~\ref{thm:r-estm-coll-time}   we get  bounds $\dot{q}_j=O(|T-t|^{-1/3})$ and $r=O(|T-t|^{2/3})$, hence  we obtain
\begin{eqnarray*}
  r \dot{\mu}' \mathcal{C}(s) = r^{5/2} \left( O(r)O(|T-t|^{-1/3}) \right)=O(r^{7/2})O(|T-t|^{-1/3})\\
  =O(|T-t|^{7/3})O(|T-t|^{-1/3})=O(|T-t|^2) \to 0.
\end{eqnarray*}

Therefore $b''$ is bounded, so by Theorem~\ref{thm:taubwin} $b' \to 0$.

\end{proof}

\section{Shadowing of perturbed orbits near normally hyperbolic manifolds}
\label{sec:shadow-isoseg}

Consider now the partial collision orbit ${\fz}(\tau)=(r,v,s,b,\theta,\mu)(\tau)$ as a solution of the system \eqref{eq:derivative of r-McGehee}-\eqref{eq:derivativeOfOmega} with separate equations for $\theta$, $\mu$.
The system has some non-autonomous terms which were exponentially estimated in Prop.\ \ref{prop:deltaEstimation}. Without these terms, Equations \eqref{eq:derivative of r-McGehee}-\eqref{eq:derivativeOfOmega} become Equations \eqref{eq:r'-auto}-\eqref{eq:zw'-auto}, which are equivalent to Eqs (3.1) from \cite{Moeckel2023}:
\[
\begin{array}{cl|l}
	& \text{Eqs \eqref{eq:r'-auto}-\eqref{eq:zw'-auto}} & \text{non-autonomous terms} \\
	\hline
	r' =& vr, &\\
    v' =& \dfrac{v^2}{2} + F(s, \zw) -V(s) &+ \dfrac{\mu^2}{r} + r^2\dfrac{\partial U_{\text{ext}}}{\partial r}, \\
    s' =& \zw, &\\
    \zw' =& -\dfrac{v\zw}{2} + \tilde{\nabla}V(s) + \frac{1}{2}\tilde{\nabla}F(s,\zw) - \mathcal{A}^{-1}(s)\left(D\mathcal{A}(s)\zw\right)\zw &+ r^2\dfrac{\partial U_{\text{ext}}}{\partial s}
    -  r\dot{\mu}\dfrac{\mathcal{A}(s)^{-1} \mathcal{B}(s)}{\|(s,1)\|^2}.
\end{array}
\]

Following \cite{Moeckel2023}, by \emph{collision manifold} we will understand the manifold $\mathcal{C}_0:=\{r=0\}$, invariant for \eqref{eq:r'-auto}-\eqref{eq:zw'-auto}. As we have seen in Subsection \ref{sec:convergence}, the partial collision orbit ${\fz}$ converges to an equilibrium point $P$ contained in $\mathcal{C}_0$, and $P$ is of the form $\{(r=0,v,s,b=0), s\in \textrm{CC}\}$, where $\textrm{CC}$ is the set of scaled central configurations.

The main goal of this section is to prove that, for a partial collision solution ${\fz}(\tau)$ of the perturbed equations \eqref{eq:derivative of r-McGehee}-\eqref{eq:derivativeOfOmega} and $\tau$ large enough, there exist an `unperturbed' solution $z_1$ of \eqref{eq:r'-auto}-\eqref{eq:zw'-auto} converging to $P$, and such that
\begin{equation*}
   \|{\fz}(\tau) - z_1(\tau)\| \leq C_0 e^{-\gamma_0 \tau}, \quad \gamma_0 > 0,
\end{equation*}

To find such a $z_1$, we will use a topological tool called \emph{isolating segment}.


\subsection{Isolating segments}

Recall from \cite{Zsegments} (originally for equations periodic in time: \cite{Srz}) the notion of isolating segment and its basic topological properties.

For a non-autonomous ODE with $x\in \mathbb{R}^n$ and a $\mathcal{C}^1$ function $f:\mathbb{R}\times\mathbb{R}^n\to\mathbb{R}^n$
\begin{equation}\label{appSegments:ode}
\dot{x}=f(t,x)
\end{equation}
we define
\[
\varphi_{(t_0,\tau)}(x_0) = x(t_0,x_0;t_0+\tau)\text{,}
\]
where by $x(t_0,x_0;\cdot)$ we understand the solution of \eqref{appSegments:ode} such that $x(t_0,x_0;t_0)=x_0$ and the domain of $\tau$ depends on $(t_0,x_0)$. Then we define the local flow for \eqref{appSegments:ode} on $\mathbb{R}\times\mathbb{R}^n$ as
\[
\Phi_{\tau} (t_0,x_0) = (t_0+\tau,\varphi_{(t_0,\tau)}(x_0)).
\]

We shall denote the projections in the space $\mathbb{R}\times \mathbb{R}^n = \mathbb{R}\times \mathbb{R}^u\times \mathbb{R}^s \ni (t,x,y)$ by $\pi_t:\mathbb{R}\times \mathbb{R}^n\to \mathbb{R}$, and $\pi_x:\mathbb{R}\times \mathbb{R}^n\to  \mathbb{R}^u$, and $\pi_y:\mathbb{R}\times \mathbb{R}^n\to  \mathbb{R}^s$. For a set $Z\subset \mathbb{R}\times \mathbb{R}^n$ and $t\in \mathbb{R}$ we also put
\[
Z_t=\{x\in\mathbb{R}^n:(t,x)\in Z\}.
\]

We are ready to introduce the notion of

\begin{definition}[Isolating segment] An \emph{isolating segment} for the equation \eqref{appSegments:ode} is a pair of sets $(W,W^-)$, $W^-\subset W\subset \mathbb{R}\times \mathbb{R}^n$ fulfilling all the following conditions:
\begin{enumerate}[(i)]
	\item $W\cap ([a,b]\times\mathbb{R}^n)$ and $W^-\cap ([a,b]\times\mathbb{R}^n)$ are compact for all $a,b\in\mathbb{R}$;
	\item for every $\sigma \in \mathbb{R}$, $x\in \partial W_{\sigma}$
there exists $\delta>0$ such that for all $t\in (0,\delta)$
$\varphi_{(\sigma,t)}(x)\not\in W_{\sigma+t}$ or
$\varphi_{(\sigma,t)}(x)\in \operatorname{int} W_{\sigma+t}$;
	\item $W^{-}= \{ (\sigma,x)\in W: \exists \delta>0\ \forall t\in (0,\delta)\ \varphi_{(\sigma,t)}(x) \not\in W_{\sigma+t} \}$,
	\\
	and the set $W^{+}:=\operatorname{cl\,} (\partial W\setminus W^{-})$ also fulfills
	\\
	$W^{+}= \{ (\sigma,x)\in W: \exists \delta>0\ \forall t\in (0,\delta)\ \varphi_{(\sigma,-t)}(x) \not\in W_{\sigma-t} \}$;
	\item there exists $\eta>0$ such that for all $x \in W^-$ there exists
$t > 0$ such that for all $\tau \in (0,t]$  $\Phi_{\tau}(0,x)\notin
W$ and $\operatorname{dist}(\Phi_t(0,x),W) > \eta$.
\end{enumerate}
In general, an isolating segment for a non-autonomous ODE is the set in the extended phase-space 
which boundaries are proper sections of the vector field.
The sets $W^-$ and $W^+$ are called the \emph{exit} and \emph{entry set} of the isolating segment, respectively. They are commonly interpreted as the subsets on the boundary of $W$ through which the trajectories of \eqref{appSegments:ode} leave or enter the segment $W$.
\end{definition}

We also define
\begin{definition}[Exit time function]
For  the  isolating segment $W$ we define \emph{the exit time function} $\tau_{W,\varphi}$ as
\begin{equation}\label{eq:exTimeMap}
\tau_{W,\varphi} : W \ni (t_0,x_0) \mapsto \sup \{ t \geq 0: \:
    \forall s \in [0,t] \: (t_0+s,\varphi_{(t_0,s)}(x_0)) \in W \} \in [0,\infty].
\end{equation}
\end{definition}
A well-known fact is that the map $\tau_{W,\varphi}$ is continuous by the Wa\.{z}ewski Retract Theorem \cite{waz}.


\subsection{General construction of a shadowing solution}

Consider a general situation of an autonomous equation for $z\in \mathbb{R}^n$ and a $\mathcal{C}^1$ map $f:\mathbb{R}^n\to\mathbb{R}^n$:
\begin{equation}
  z'=f(z),  \label{eq:z'=f}
\end{equation}
with $f(0)=0$. Let $\varphi(t,z)$ be the local dynamical system induced on $\mathbb{R}^n$ by (\ref{eq:z'=f}).

After performing the linearization of $f$ in $0$, we obtain some coordinates $(x,y)\in\mathbb{R}^u\times\mathbb{R}^s = \mathbb{R}^n$ such that (\ref{eq:z'=f}) becomes
\begin{eqnarray}
  x'&=&f_x(x,y)=A_{cu} x + N_x(x,y)  \label{eq:x'cdir}  \\
  y'&=&f_y(x,y)=A_{s} y + N_y(x,y), \label{eq:y'cdir}
\end{eqnarray}
with
\begin{equation}
  N(x,y)=O(\|(x,y)\|^2), \quad DN(x,y)=O(\|(x,y)\|), \label{eq:NO()}
\end{equation}
where $\mathbb{R}^s\ni y$ is the stable direction for linearization  and $\mathbb{R}^u\ni x$ represents the center-unstable subspace.
Therefore, we can assume that (we use the same notation `$\cdot$' for the standard inner products in $\mathbb{R}^s$, $\mathbb{R}^u$ and $\mathbb{R}^n$):
\begin{eqnarray}
   \mu_\text{log}(A_s):=\sup_{\|y\|=1} A_sy\cdot y < 0, \qquad m_l(A_{cu}):=\inf_{\|x\|=1} A_{cu}x\cdot x\geq 0.  \label{eq:logN-lin}
\end{eqnarray}

Let $U=\overline{B}_n(0,R)$ be the closed ball of radius $R$ in $\mathbb{R}^n$ (shortly: $\overline{B}_n$ for unit ball).
We define (the names of constants are suggested by \cite{CZ15}, we set $L=1$)
\begin{eqnarray*}
 \overrightarrow{\mu}= \overrightarrow{\mu _{s,1}}&=&\sup_{z\in U}\left\{ \mu_\text{log}\left( \frac{\partial
f_{y}}{\partial y}(z)\right) +
\left\Vert \frac{\partial f_{y}}{\partial x}(z)\right\Vert \right\}, \\
\overrightarrow{\xi}=\overrightarrow{\xi _{cu,1,U}}&=&\inf_{z\in U} m_{l}\left( \frac{\partial
f_{x}}{\partial x}(z)\right) -
\sup_{z\in U}\left\Vert \frac{\partial f_{x}}{\partial y}(z)\right\Vert\text{.}
\end{eqnarray*}

From (\ref{eq:x'cdir}, \ref{eq:y'cdir}, \ref{eq:NO()}) and \eqref{eq:logN-lin} it follows that for sufficiently small $R$ holds
\begin{equation}
   \overrightarrow{\mu} < 0, \quad  \overrightarrow{\mu} < \overrightarrow{\xi}.  \label{eq:cm-cone-cond}
\end{equation}
What is more, $\overrightarrow{\xi}$ is negative and $\overrightarrow{\xi}\to 0$ with $R\to 0$.

\begin{theorem}\label{theo:shadowing theorem}
Assume that $R$ is such that \eqref{eq:cm-cone-cond} holds.
Consider a perturbed orbit of \eqref{eq:z'=f} fulfilling
\begin{equation*}
  z_p'(t)=f(z_p(t)) + \delta(t),
\end{equation*}
such that $z_p(t)=(x_p(t),y_p(t)) \in U$, and for $t \in [0,\infty)$:
\begin{equation*}
  \|\delta(t)\| \leq a e^{\alpha t}
\end{equation*}
with some constants $a$ and $\alpha <0$, $\alpha <  \overrightarrow{\xi}$.

Then for any $\gamma <0$, such that
\begin{equation}
  \overrightarrow{\mu} < \gamma < \overrightarrow{\xi},  \quad \alpha < \gamma,  \label{eq:cm-ode-gamma}
\end{equation}
and any $r>0$
 there exists $z_1 \in U$ and $t_0\geq 0$, such that $\varphi(t,z_1) \in U$ for $t \in [0,\infty)$ and
\begin{equation*}
  \|z_p(t_0+t) - \varphi(t,z_1)\| \leq re^{\gamma t}, \quad t \geq t_0.
\end{equation*}
\end{theorem}

\begin{proof}
Let us fix $\gamma<0$ satisfying \eqref{eq:cm-ode-gamma} and any $r>0$. We will construct an isolating segment $W=\{W_t\}_{t \in \mathbb{R}_+}$ on extended phase-space $\mathbb{R}\times \mathbb{R}^u\times \mathbb{R}^s\ni(t,x,y)$. 
\begin{eqnarray*}
 W_t&=&z_p(t) + \overline{B}_u(0,re^{\gamma t}) \times \overline{B}_s(0,re^{\gamma t}), \\
 W_t^-&=&z_p(t) + \partial(\overline{B}_u(0,re^{\gamma t})) \times \overline{B}_s(0,re^{\gamma t}), \\
 W_t^+&=&z_p(t) + \overline{B}_u(0,re^{\gamma t}) \times \partial \overline{B}_s(0,re^{\gamma t}).
\end{eqnarray*}
We need to show that on $W^-$ the vector field $f$ points outside of $W$ and on $W^+$ it points to the interior of $W$.

Observe that for some $t_U \geq 0$, $W_{t} \subset U$ for $t\geq t_U$. From now on, we will assume $t\geq t_U$.

Let us set
\begin{eqnarray*}
  L^-(t,z):=\| x - x_p(t) \|^2 - r^2e^{2\gamma t}, \\
  L^+(t,z):=\|y - y_p(t)\|^2 - r^2e^{2\gamma t}.
\end{eqnarray*}
It is easy to see for $W^-$ to be  the exit set we need
\begin{equation}
  DL^-(t,z) \cdot (1,f) > 0, \quad z \in W_t^-, \label{eq:dL-}
\end{equation}
and analogously for $W^+$ to be an entry set we need that
\begin{equation}
  DL^+(t,z) \cdot (1,f) < 0, \quad z \in W_t^+. \label{eq:dL+}
\end{equation}

We start with condition (\ref{eq:dL-}). We have
\begin{eqnarray*}
  \frac{\partial L^-}{\partial t}(t,z)&=& 2(x-x_p(t))\cdot(-x'_p(t)) - 2\gamma r^2e^{2\gamma t}\\
   &=& -2(x-x_p(t))\cdot (f_x(z_p(t)) + \delta_x(t))  - 2\gamma r^2e^{2\gamma t}, \\
  \frac{\partial L^-}{\partial x}(t,z) h &=& 2(x-x_p(t))\cdot h, \\
   \frac{\partial L^-}{\partial y}(t,z) h &=& 0,
\end{eqnarray*}
therefore for $z\in W_t^-$  we have
\begin{multline*}
    DL^-(t,z) \cdot (1,f) = -2(x-x_p(t))\cdot(f_x(z_p(t))+\delta_x(t))  - 2\gamma r^2e^{2\gamma t}  + 2(x-x_p(t))\cdot f_x(z) =
    \\
      = 2(x-x_p(t))\cdot (f_x(z)-f_x(z_p(t))) -2(x-x_p(t))\cdot\delta_x(t)  - 2\gamma r^2e^{2\gamma t} \geq
     \\
      \geq 2\inf_{z \in W_t^-} (x-x_p(t))\cdot (f_x(z)-f_x(z_p(t)))  - 2re ^{\gamma t} \|\delta(t)\| - 2\gamma r^2 e^{2\gamma t}.
\end{multline*}
Let us set $z-z_p(t)=\bar{z}$, then $\bar{z} \in  \partial(\overline{B}_u(0,re^{\gamma t})) \times \overline{B}_s(0,re^{\gamma t})$
and
\begin{equation*}
   \inf_{z \in W_t^-} (x-x_p(t))\cdot (f_x(z)-f_x(z_p(t))) = \inf_{\substack{
   \|\bar{x}\|=r e^{\gamma t}\\ \|\bar{y}\| \leq r e^{\gamma t}}}\bar{x}\cdot(f_x(z_p(t)+\bar{z})-f_x(z_p(t)))
\end{equation*}
since
\[
f_x(z_p(t)+\bar{z})-f_x(z_p(t)) = \bar{x} \int_0^1 \frac{\partial f_x}{\partial x}(z_p(t)+s\bar{z}) ds
  + \bar{y}\int_0^1 \frac{\partial f_x}{\partial y}(z_p(t)+s\bar{z}) ds ,
\]
hence
\[
  \bar{x}\cdot(f_x(z_p(t)+\bar{z})-f_x(z_p(t))) \geq \left(m_l\left(\frac{\partial f_x}{\partial x}(U) \right)   - \left\|\frac{\partial f_x}{\partial y}(U) \right\| \right)   r^2 e^{2\gamma t}
  \geq \overrightarrow{\xi}  r^2 e^{2\gamma t}.
\]
Therefore we obtain for $z\in W_t^-$
\begin{multline*}
   DL^-(t,z) \cdot (1,f) \geq 2 \overrightarrow{\xi}  r^2 e^{2\gamma t} - 2r a e^{(\gamma + \alpha) t} - 2\gamma r^2 e^{2\gamma t}=
   \\
     = 2 r^2 (\overrightarrow{\xi} - \gamma) e^{2\gamma t} - 2r a e^{\gamma t} e^{\alpha t}= 2r e^{\gamma t} \left(r (\overrightarrow{\xi} - \gamma)e^{\gamma t} - a e^{\alpha t} \right),
\end{multline*}
and we see that $ DL^-(t,z) \cdot (1,f)>0$ for sufficiently large $t > t^1_0$.

Now we check condition (\ref{eq:dL+}).
 We have
\begin{eqnarray*}
  \frac{\partial L^+}{\partial t}(t,z)&=& 2(y-y_p(t))\cdot(-y'_p(t)) - 2\gamma r^2e^{2\gamma t}\\
   &=& -2(y-y_p(t))\cdot f_y(z_p(t)) + \delta_y(t))  - 2\gamma r^2e^{2\gamma t}, \\
  \frac{\partial L^+}{\partial x}(t,z) h &=& 0, \\
   \frac{\partial L^+}{\partial y}(t,z) h &=& (y-y_p(t))\cdot h,
\end{eqnarray*}
Therefore for $z\in W_t^+$  we have
\begin{multline*}
    DL^+(t,z) \cdot (1,f) = -2(y-y_p(t))\cdot (f_y(z_p(t))+\delta_y(t))  - 2\gamma r^2e^{2\gamma t}   + 2(y-y_p(t))\cdot f_y(z) =
    \\
      = 2(y-y_p(t))\cdot(f_y(z)-f_y(z_p(t))) -2(y-y_p(t))\cdot \delta_y(t)  - 2\gamma r^2e^{2\gamma t} \leq
      \\
      \leq 2\sup_{z \in W_t^+} (y-y_p(t))\cdot(f_y(z)-f_y(z_p(t)))  + 2re ^{\gamma t} \|\delta(t)\| - 2\gamma r^2 e^{2\gamma t}.
\end{multline*}
Let us now set $z-z_p(t)=\bar{z}$, then $\bar{z} \in  \overline{B}_u(0,re^{\gamma t}) \times \partial \overline{B}_s(0,re^{\gamma t})$
and
\begin{equation*}
   \sup_{z \in W_t^+} (y-y_p(t))\cdot(f_y(z)-f_y(z_p(t))) = \sup_{\substack{
   \|\bar{x}\|\leq r e^{\gamma t}\\ \|\bar{y}\| = r e^{\gamma t}}}\bar{y}\cdot(f_y(z_p(t)+\bar{z})-f_y(z_p(t)))
\end{equation*}
since
\[
f_y(z_p(t)+\bar{z})-f_y(z_p(t))=  \bar{x}\int_0^1 \frac{\partial f_y}{\partial x}(z_p(t)+s\bar{z}) ds
  + \bar{y}\int_0^1 \frac{\partial f_y}{\partial y}(z_p(t)+s\bar{z}) ds ,
\]
hence
\[
  \bar{y}\cdot(f_y(z_p(t)+\bar{z})-f_y(z_p(t))) \leq \left(\mu_\text{log}\left(\frac{\partial f_y}{\partial y}(U) \right)   + \left\|\frac{\partial f_y}{\partial x}(U) \right\| \right)   r^2 e^{2\gamma t}
  \leq \overrightarrow{\mu}  r^2 e^{2\gamma t}.
\]
Therefore we obtain for $z\in W_t^+$
\begin{multline*}
   DL^+(t,z) \cdot (1,f) \leq  2 \overrightarrow{\mu}  r^2 e^{2\gamma t} + 2r a e^{(\gamma + \alpha) t} - 2\gamma r^2 e^{2\gamma t}=
   \\
     = 2 r^2 (\overrightarrow{\mu} - \gamma) e^{2\gamma t} + 2r a e^{\gamma t} e^{\alpha t}= -2r e^{\gamma t} \left(r (\gamma - \overrightarrow{\mu} )e^{\gamma t} - a e^{\alpha t} \right),
\end{multline*}
and we see that $ DL^+(t,z) \cdot (1,f)<0$ for sufficiently large $t > t^2_0$.

This shows that $W$ is an isolating segment for $ t \geq t_0:=\max\{t^1_0,t^2_0,t_U\}$.

Now we follow the standard topological argument used in proving the existence of an orbit contained in the isolating segment.


\begin{lemma}\label{lem:retraction}
There exists $z_1 \in W_{t_0}$, such that $\varphi(t,z_1) \in W_{t_0+t}$, for $ t \in [t_0,+\infty)$.
\end{lemma}

\begin{proof}
First, observe that for $u=0$ (i.e. $W^-=\emptyset$) no orbit entering the segment $W$ through $W_{t_0}$ ever exits it, so every $z\in W_{t_0}$ fulfills the thesis of the lemma.

Assume now $u>0$ and suppose to the contrary that for every $z \in W_{t_0} = z_p(t_0) + \overline{B}_u(0,re^{\gamma t_0}) \times \overline{B}_s(0,re^{\gamma t_0})$ the forward orbit  $\varphi([t_0,\infty),z)$ exits the segment $W$, that is, its exit time function \eqref{eq:exTimeMap}
\[
\tau(z):=\tau_{W,\varphi}(t_0,z) = \sup \{ t \geq 0: \:
    \forall s \in [0,t] \: (t_0+s,\varphi_{(t_0,s)}(z)) \in W \} <\infty.
\]

 Under the above assumption we will construct a deformation retraction of  $W_{t_0}$ onto  $W_{t_0}^-$, this will lead to a contradiction.

Consider a function
\[
h:[0,1]\times W_{t_0}\ni (s,z) \mapsto e^{-\gamma s\tau(z)}[\varphi(s\tau(z),z)-z_p(t_0+s\tau(z))]+z_p(t_0).
\]
From the properties of isolating segments, $h$ is well-defined and continuous. We claim that:

\begin{enumerate}
	\item $h([0,1]\times W_{t_0})\subset W_{t_0}$.
	\\
	Indeed, for a fixed $(s,z)\in [0,1]\times W_{t_0}$:
	\begin{multline*}
	\|\pi_x(h(s,z)-z_p(t_0))\| =e^{-\gamma s\tau(z)}\|\pi_x([\varphi(s\tau(z),z)-z_p(t_0+s\tau(z))])\|\leq
	\\
	\leq e^{-\gamma s \tau(z)}re^{\gamma (t_0+s\tau(z))}=re^{\gamma t_0},
	\end{multline*}
	because $\varphi(s\tau(z),z)\in W_{t_0+s\tau(z)} = z_p(t_0+s\tau(z)) + \overline{B}_u(0,re^{\gamma(t_0+s\tau(z))}) \times \overline{B}_s(0,re^{\gamma(t_0+s\tau(z))})$. Analogously we see that
	\[
	\|\pi_y(h(s,z)-z_p(t_0))\| \leq re^{\gamma t_0}.
	\]
	\item $h(0,z)= e^{0}[\varphi(0,z)-z_p(t_0)]+z_p(t_0)=z$,\quad so $h(0,\cdot)$ is the identity function on $W_{t_0}$.
	\item\label{3} For $z\in W^-_{t_0}$ and any $s\in[0,1]$, $h(s,z)= z$.
	\\	
	This is obvious because for such $z$, $\tau(z)=0$.
	\item $h(\{1\}\times W_{t_0})= W^-_{t_0}$.
	\\
	Regarding \ref{3}, it is sufficient to show that for any $z\in W_{t_0}$, $h(1,z)\in W^-_{t_0}$:
	\begin{multline*}
	\|\pi_x(h(1,z)-z_p(t_0))\| =e^{-\gamma \tau(z)}\|\pi_x([\varphi(\tau(z),z)-z_p(t_0+\tau(z))])\|=
	\\
	= e^{-\gamma \tau(z)}re^{\gamma (t_0+\tau(z))}=re^{\gamma t_0},
	\end{multline*}
	because $\varphi(\tau(z),z)\in W^-_{t_0+\tau(z)} = z_p(t_0+\tau(z)) + \partial(\overline{B}_u(0,re^{\gamma(t_0+\tau(z))})) \times \overline{B}_s(0,re^{\gamma(t_0+\tau(z))})$.
\end{enumerate}

From observations 1--4 we conclude that $h$ would be a deformation retraction of $W_{t_0}$ onto $W^-_{t_0}$, which is impossible because the sets have different reduced singular homology groups, in particular:
\begin{equation*}
\bar H_{u-1}(W_{t_0})=\bar H_{u-1}(\overline{B}_u\times\overline{B}_s) = \{0\}
\text{, \qquad }
\bar H_{u-1}(W^-_{t_0})=\bar H_{u-1}(\partial\overline{B}_u\times\overline{B}_s) = \mathbb{Z}. 
\end{equation*}
We obtain a contradiction which proves that there must exist a point $z_1\in W_{t_0}$ for which $\tau(z_1)=+\infty$. Therefore, its forward orbit must be contained in the isolating segment $W$.
\end{proof}

The segment was constructed in such a way that $z_1$ from Lemma \ref{lem:retraction} fulfills the thesis of Theorem~\ref{theo:shadowing theorem}.
\end{proof}

\subsection{Shadowing partial collision orbit by an orbit of autonomous system   \eqref{eq:r'-auto}-\eqref{eq:zw'-auto}}

Let us go back to our particular case of partial collision orbit $\fz$. We claim that the systems \eqref{eq:derivative of r-McGehee}-\eqref{eq:derivativeOfOmega} and \eqref{eq:r'-auto}-\eqref{eq:zw'-auto} fulfill the assumptions of Theorem~\ref{theo:shadowing theorem}:
\begin{itemize}
	\item The system \eqref{eq:r'-auto}-\eqref{eq:zw'-auto} is linearized around the limit point $P$ (see \cite[Sec.\ 4]{Moeckel2023} for analysis of the center manifold) and we define the constants $ \overrightarrow{\mu}$, $\overrightarrow{\xi}$ fulfilling \eqref{eq:cm-cone-cond} in some closed ball $U=\overline{B}_{4k-6}(P,R)$ of sufficiently small radius $R$;
	\item We regard the system \eqref{eq:derivative of r-McGehee}-\eqref{eq:derivativeOfOmega} as a non-autonomous perturbation of the system \eqref{eq:r'-auto}-\eqref{eq:zw'-auto} with the perturbation exponential estimation $\|\delta(t)\|\leq C e^{-Et}$ given by Proposition \ref{prop:deltaEstimation};
	\item Recall that the constant 	
	  $\overrightarrow{\xi} =\inf\limits_{z\in U} \inf\limits_{\|x\|=1} \left( \frac{\partial f_{x}}{\partial x}(z) x\right)\cdot x -
\sup\limits_{z\in U}\left\Vert \frac{\partial f_{x}}{\partial y}(z)\right\Vert$, where $f$ is the rhs of \eqref{eq:r'-auto}-\eqref{eq:zw'-auto} and by $x$ and $y$ we understand the local coordinates on the center-unstable and stable manifolds, respectively.
As we mentioned at the properties \eqref{eq:cm-cone-cond} of $\overrightarrow{\xi}$, we can make it as close to zero as we wish, decreasing $R$. Hence, we can assume
$-E<\overrightarrow{\xi}$.
\end{itemize}

Therefore, from Theorem \ref{theo:shadowing theorem}, we may state that

\begin{corollary}\label{cor:shadowing}
Let ${\fz}(\tau)$, $\tau\in[0,\infty)$, be a partial collision solution of the system \eqref{eq:derivative of r-McGehee}-\eqref{eq:derivativeOfOmega}. Then there exist a solution $z_1(\tau)$ of the system \eqref{eq:r'-auto}-\eqref{eq:zw'-auto} and $C_0, \gamma_0 > 0$ such that for $\tau$ large enough
\begin{equation*}
   \|{\fz}(\tau) - z_1(\tau)\| \leq C_0 e^{-\gamma_0 \tau}.
\end{equation*}
\end{corollary}

\section{No infinite spin}
\label{sec:no-inf-spin}

Recall now some facts from \cite{Moeckel2023}, regarding the system \eqref{eq:r'-auto}-\eqref{eq:zw'-auto}.

Denote by $z_1(\tau)$ any solution of \eqref{eq:r'-auto}-\eqref{eq:zw'-auto} converging to a point $P\in \mathcal{C}_0$.
Therefore, we can assume that $z_1(\tau)$ is confined in a neighborhood $\mathcal{U}$ of $P$ and $z_1(\tau)$ is in the local center-stable manifold $W_{\mathcal{U}}^{cs}$ (see \cite{Wiggins}), 
for a large enough $\tau$.
Then, from the proof of \cite[Thm.\ 4.6]{Moeckel2023} there exists a solution $z_s(\tau) = (r_s, v_s, s_s,\zw_s)(\tau)$ contained in the center manifold on the same fiber such that $\| z_1 - z_s\| \leq c_0e^{-|c_1|\tau}$ and $z_s\to P$ as $\tau\to \infty$. Finally, it is proved that the arclength of $z_s$ is finite, and consequently, the spin angle $\theta_s$ converges, i.e.,
\begin{equation}\label{eq:finite_angle_theta_s}
 \int_0^\infty |\theta'_s| d\tau = \int_0^\infty\frac{|\Omega(s_s(\tau),\zw_s(\tau))|}{\left\|(s_s,1)\right\|^2} d\tau < \infty.
\end{equation}

Consider again the partial collision solution ${\fz} = (r, v, s, b, \theta,\mu)$ of \eqref{eq:derivative of r-McGehee}-\eqref{eq:derivativeOfOmega}.
As observed by \cite{Moeckel2023}, to prove that there is no infinite spin for the partial collision orbit it is enough to ensure that
\begin{equation}\label{eq:spin_angle}
  \int_{\tau_0}^\infty |\theta'(\tau)| d\tau <\infty,
\end{equation}
this means that the angle is converging.

\label{Happy end}
Thanks to the Corollary \ref{cor:shadowing}, we are able to approximate exponentially the partial collision solution ${\fz}(\tau)$ of \eqref{eq:derivative of r-McGehee}-\eqref{eq:derivativeOfOmega} by a solution $z_1(\tau)$ of the system \eqref{eq:r'-auto}-\eqref{eq:zw'-auto} converging to some point $P$ on the collision manifold $\mathcal{C}_0$. From \cite[Thm.\ 4.6]{Moeckel2023} we conclude that the solution $z_1(\tau)$  can be approximated exponentially by a solution $z_s(\tau)$ of \eqref{eq:r'-auto}-\eqref{eq:zw'-auto} contained in $\mathcal{C}_0$.
As the distance between $z_1$ and $z_s$ converges to zero exponentially, the same will happen for the distance between ${\fz}$ and $z_s$, i.e.,
\begin{equation}\label{eq:sync_on_center manifold}
   \|{\fz}(\tau) - z_s(\tau)\| \leq C_1 e^{-\gamma_1 \tau}, \quad \gamma_1 > 0.
\end{equation}

Returning back to \eqref{eq:spin_angle}, from (\ref{eq:theta'}), we have
\begin{equation*}
 \theta' = \frac{\mu}{\sqrt{r}} -\frac{\Omega(s,\zw)}{\left\|(s,1)\right\|^2}.
\end{equation*}
Observe that from Lemma~\ref{lem:mu2/r}, we obtain
\begin{equation*}
  \left|  \frac{\mu}{\sqrt{r}} \right| = O(e^{-2E_1\tau}).
\end{equation*}
Hence, for $\tau$ large enough, there is a positive constant $C_2$ such that
\begin{align*}
    \int_{\tau_0}^{\infty}|\theta'(\tau)|d\tau \leq C_2\int_{\tau_0}^{\infty}e^{-2 E_1\tau }d\tau + \int_{\tau_0}^{\infty}\frac{|\Omega(s,\zw)|}{\left\|(s,1)\right\|^2}d\tau.
\end{align*}

By applying \eqref{eq:finite_angle_theta_s} and \eqref{eq:sync_on_center manifold},  we have
\begin{eqnarray*}
 \int_{\tau_0}^\infty |\theta'| d \tau  &\leq& C_2  \int_{\tau_0}^\infty e^{-2 E_1\tau } d\tau + \int_{\tau_0}^\infty\frac{|\Omega(s_s(\tau),\zw_s(\tau))|}{\left\|(s_s,1)\right\|^2} d\tau \\
   & & + \int_{\tau_0}^\infty \sup\limits_{z \in \mathcal{U}_1} \left\| D \left(\frac{\Omega(s,\zw)}{\left\|(s,1)\right\|^2}\right) \right\| |z_s(\tau) - {\fz}(\tau)|d\tau  \\
   &\leq&  C_2  \int_{\tau_0}^\infty e^{-2 E_1\tau } d\tau  + \int_{\tau_0}^\infty\frac{|\Omega(s_s(\tau),\zw_s(\tau))|}{\left\|(s,1)\right\|^2} d\tau \\
 & &  +C_1\sup\limits_{z \in \mathcal{U}_1} \left\| D\left(\frac{\Omega(s,\zw)}{\left\|(s,1)\right\|^2}\right)\right\| \int_{\tau_0}^\infty e^{-\gamma_1 \tau}d\tau < \infty,
\end{eqnarray*}
where $\mathcal{U}_1$ is some bounded convex set containing both ${\fz}(\tau)$ and $z_s(\tau)$ for $\tau$ large enough.
Therefore, we have proved the following theorem.

\begin{theorem}\label{th:main}
    A partial collision orbit on the planar $n$-body problem that converges to an isolated CC does not have an infinite spin.
\end{theorem}


\appendix

\section{Computation of $\dot{\omega}$}
\label{sec:dotomega}
The goal of this section is to compute $\dot{\omega}$ for Lagrangian (\ref{eq:LagMMcoord}).
As explained in Subsec. \ref{secsub:eq-motion}, the equation for $\dot{\omega}$ is derived from the Lagrange equation:
\begin{align*}
    \frac{d}{dt}\frac{\partial \mathcal{L} }{\partial \omega}- \frac{\partial \mathcal{L} }{\partial s} = 0.
\end{align*}
To make sense of various formulas, observe first that $\frac{\partial \mathcal{L}}{\partial s}$ or  $\frac{\partial \mathcal{L}}{\partial \omega}$ are both forms and should therefore be expressed as row-vectors.

The terms containing $\dot{\theta}$ and are involved in $\mu$ can be gathered below in $K$ given by
\begin{eqnarray*}
  K(r,s,\omega,\dot{\theta})=\frac{r^2 \Omega(s,\omega)^2}{2\left\|(s,1)\right\|^4} + \frac{r^2\dot{\theta}\Omega(s,\omega)}{\left\|(s,1)\right\|^2}=\frac{1}{\|(s,1)\|^2}\left( \frac{r^2 \Omega(s,\omega)^2}{2\left\|(s,1)\right\|^2}  + r^2\dot{\theta}\Omega(s,\omega) \right).
\end{eqnarray*}

\subsection{Computation of  $\frac{\partial \mathcal{L}}{\partial \omega}$}
Let us begin the calculation of  $\frac{\partial \mathcal{L}}{\partial \omega}$.
The term $\frac{r^2}{2} F(s,\omega)$ is easy and we obtain
\begin{eqnarray*}
  \frac{\partial }{\partial \omega} \left( \frac{r^2}{2} F(s,\omega) \right)=r^2 \omega^T \mathcal{A}(s),
\end{eqnarray*}
and for $K(r,s,\omega,\dot{\theta})$, using \eqref{eq:mu-dot-theta}, we obtain
\begin{eqnarray*}
   \frac{\partial }{\partial \omega}K=
   \frac{1}{\|(s,1)\|^2}\left( \frac{r^2 \Omega(s,\omega)}{\left\|(s,1)\right\|^2}  + r^2\dot{\theta} \right)\frac{\partial \Omega(s,\omega)}{\partial \omega}
   = \mu \frac{\mathcal{B}(s)^T}{\|(s,1)\|^2}.
\end{eqnarray*}
Therefore,
\begin{eqnarray}\label{eq:partial_omega_L}
  \frac{\partial \mathcal{L}}{\partial \omega}= r^2 \omega^T \mathcal{A}(s) + \mu \frac{\mathcal{B}(s)^T}{\|(s,1)\|^2}.
\end{eqnarray}

\subsection{Computation of $\frac{\partial \mathcal{L}}{\partial s}$}

The most demanding or tricky part is coming from $K$. We have
\begin{eqnarray*}
    \frac{\partial }{\partial s}K= D_s \left(\frac{1}{\|(s,1)\|^2} \right) \left( \frac{r^2 \Omega(s,\omega)^2}{2\left\|(s,1)\right\|^2}  + r^2\dot{\theta} \Omega(s,\omega)\right) + \frac{1}{\|(s,1)\|^2}\left( \frac{r^2 \Omega(s,\omega)}{\left\|(s,1)\right\|^2}  + r^2\dot{\theta} \right)\frac{\partial \Omega(s,\omega)}{\partial s} \\
    + \frac{1}{\|(s,1)\|^2} \frac{r^2 \Omega(s,\omega)^2}{2} D_s\left(\frac{1}{\|(s,1)\|^2} \right) \\
   = D_s\left(\frac{1}{\|(s,1)\|^2} \right)\left(  \frac{r^2 \Omega(s,\omega)^2}{\left\|(s,1)\right\|^2}  + r^2\dot{\theta} \Omega(s,\omega) \right) + \mu \frac{\mathcal{B}'(s)^T \omega}{\|(s,1)\|^2} \\
   = -2\frac{D_s \|(s,1)\|}{\|(s,1)\|^3}\left(\frac{r^2 \Omega(s,\omega)^2}{\left\|(s,1)\right\|^2} + \left(\mu - \frac{r^2 \Omega(s,\omega)}{\left\|(s,1)\right\|^2}\right)\Omega(s,\omega) \right) + \mu \frac{\mathcal{B}'(s)^T \omega}{\|(s,1)\|^2} \\
   = -2\frac{D_s \|(s,1)\|}{\|(s,1)\|^3} \mu \Omega(s,\omega) + \mu \left( \frac{\mathcal{B}'(s)^T \omega}{\|(s,1)\|^2} \right) \\
   = \mu \frac{\partial }{\partial s}\left(\frac{\Omega(s,\omega)}{\|(s,1)\|^2} \right).
\end{eqnarray*}

The other terms in $\mathcal{L}$ can be easily obtained. Therefore,
\begin{eqnarray*}
  \frac{\partial \mathcal{L}}{\partial s} = \frac{r^2}{2} D_s F(s,\omega) + \mu D_s\left(\frac{\Omega(s,\omega)}{\|(s,1)\|^2} \right) + \frac{D_s V(s)}{r} + D_s U_{\text{ext}}.
\end{eqnarray*}

\subsection{Computation of $\frac{d}{dt}\frac{\partial \mathcal{L}}{\partial \omega}$ and equation for $\dot{\omega}$}

We are just one step away from obtaining the equation for $\dot{\omega}$. Now, we calculate the derivative of \eqref{eq:partial_omega_L} with respect to time $t$, i.e.,
\begin{eqnarray*}
  \frac{d}{dt}\frac{\partial \mathcal{L}}{\partial \omega}^T =2 r \dot{r} \mathcal{A}(s) \omega + r^2 (\nabla_s \mathcal{A}(s) \dot{s}) \omega + r^2 \mathcal{A}(s) \dot{\omega} + \dot{\mu} \frac{\mathcal{B}(s)}{\|(s,1)\|^2} + \mu \nabla_s \left(\frac{\mathcal{B}(s)}{\|(s,1)\|^2}\right)\dot{s} \\
  = 2 r \rho \mathcal{A}(s) \omega + r^2 (\nabla_s \mathcal{A}(s) \omega) \omega + r^2 \mathcal{A}(s) \dot{\omega} + \dot{\mu} \frac{\mathcal{B}(s)}{\|(s,1)\|^2} + \mu \nabla_s \left(\frac{\mathcal{B}(s)}{\|(s,1)\|^2}\right)\omega.
\end{eqnarray*}

Finally, we derive the following equation for $\dot{\omega}$:
\begin{eqnarray*}
  \dot{\omega} = -\frac{2 \rho}{r}  \omega - \mathcal{A}(s)^{-1}(\nabla_s \mathcal{A}(s) \omega) \omega - \frac{\mu}{r^2} \mathcal{A}(s)^{-1}\nabla_s \left(\frac{\mathcal{B}(s)}{\|(s,1)\|^2}\right)\omega
    - \frac{\dot{\mu}}{r^2}\mathcal{A}(s)^{-1}\frac{\mathcal{B}(s)}{\|(s,1)\|^2} \\
    + \frac{1}{2}\mathcal{A}(s)^{-1}\nabla_s F(s,\omega) + \frac{1}{r^3} \mathcal{A}(s)^{-1}\nabla_s V(s) + \frac{1}{r^2} \mathcal{A}(s)^{-1} \nabla_s U_{\text{ext}}+ \frac{\mu}{r^2} \mathcal{A}(s)^{-1}  \nabla_s\left(\frac{\Omega(s,\omega)}{\|(s,1)\|^2} \right).
\end{eqnarray*}
By using the notation $\widetilde{\nabla}_s = \mathcal{A}(s)^{-1} \nabla_s$ and noting $\nabla_s\left(\frac{\Omega(s,\omega)}{\|(s,1)\|^2}\right) =\nabla_s \left(\frac{\mathcal{B}(s)}{\|(s,1)\|^2}\right)\omega$,  we can rewrite the aforementioned expression as follows:
\begin{eqnarray*}
 \dot{\omega} = -\frac{2 \rho}{r}  \omega - (\widetilde{\nabla}_s \mathcal{A}(s) \omega) \omega
    + \frac{1}{2}\widetilde{\nabla}_s F(s,\omega) + \frac{1}{r^3} \widetilde{\nabla}_s V(s) + \frac{1}{r^2} \widetilde{\nabla}_s U_{\text{ext}}\\
    + \frac{\mu}{r^2} \mathcal{A}(s)^{-1}  \nabla_s\left(\frac{\Omega(s,\omega)}{\|(s,1)\|^2} \right) - \frac{\mu}{r^2} \mathcal{A}(s)^{-1}\nabla_s \left(\frac{\mathcal{B}(s)}{\|(s,1)\|^2}\right)\omega
    - \frac{\dot{\mu}}{r^2}\mathcal{A}(s)^{-1}\frac{\mathcal{B}(s)}{\|(s,1)\|^2} \\
    = -\frac{2 \rho}{r}  \omega - (\widetilde{\nabla}_s \mathcal{A}(s) \omega) \omega
    + \frac{1}{2}\widetilde{\nabla}_s F(s,\omega) + \frac{1}{r^3} \widetilde{\nabla}_s V(s) + \frac{1}{r^2} \widetilde{\nabla}_s U_{\text{ext}}  - \frac{\dot{\mu}}{r^2}\mathcal{A}(s)^{-1}\frac{\mathcal{B}(s)}{\|(s,1)\|^2}.
\end{eqnarray*}

\section{Asymptotic bounds for partial collisions}
\label{sec:est-coll}

In this section, we do not assume a planar case; the results derived here apply to both planar and spatial cases. Our objective is to establish asymptotic bounds
for certain physical quantities of a system within a cluster experiencing a partial collision.

Consider $G_1=\mathcal{G}, G_2,\dots, G_{N_c}$, clusters for (possibly multiple) partial collisions, with $\mathcal{G}$ containing at least two bodies.
 It is important to note that trajectory $q_i(t)$ converges to $L_i$ as $t \nearrow T$, and $L_i=L_j$ if and only if indices $i$ and $j$ belong to the same cluster.
 For a cluster $G_k$, we will often use $L_{G_k}=L_j$ for $j \in G_k$.


Let $h$ be the energy of the colliding solution. Let $K=K(t)$ be the kinetic energy, then $h=K(t)-U(t)$. Recall that $M_\mathcal{G}$ is the total mass of bodies in cluster $\mathcal{G}$. By $M_\text{tot}=\sum_{i=1}^n m_i$, we denote the total mass of all $n$ bodies.
For a given subset $\mathcal{S} \subset \{1,\dots,n\}$, we set
\begin{eqnarray*}
K_\mathcal{S}(t)&=&\sum_{i \in \mathcal{S}} \frac{m_i \dot{q}_i(t)^2}{2}, \\
U_{\mathcal{S}}(t)&=&\sum_{j<k, j,k \in \mathcal{S}} U_{j,k}(q_j(t),q_k(t)), \\
H_\mathcal{S}(t)&=&K_\mathcal{S}(t)- U_\mathcal{S}(t), \\
M_\mathcal{S}&=&\sum_{i \in \mathcal{S}} m_i, \\
c_\mathcal{S}(t)&=&M_\mathcal{S}^{-1} \sum_{i \in \mathcal{S}} m_i q_i(t), \\
\mu_\mathcal{S}&=& \sum_{i \in \mathcal{S}} m_i q_i(t) \times \dot{q}_i(t), \\
\mu^0_{\mathcal{S}}&=& \sum_{i \in \mathcal{S}} m_i \left(q_i(t)-c_\mathcal{S}(t) \right)\times \left(\dot{q}_i(t) - \dot{c}_\mathcal{S}(t)\right), \\
I^0_\mathcal{S}&=&\sum_{i\in \mathcal{S}} m_i |q_i-c_\mathcal{S}|^2, \\
r_\mathcal{S}&=&\sqrt{I^0_\mathcal{S}}.
\end{eqnarray*}
Here, $H_\mathcal{S}$ represents the ``internal'' energy of bodies in $\mathcal{S}$, and $c_\mathcal{S}$  is the center of mass of bodies in $\mathcal{S}$. It is worth noting that if $\mathcal{S}=\{1,\dots,n\}$, then $H_\mathcal{S}(t)=h$, $K_\mathcal{S}(t)=K(t)$,
$U_\mathcal{S}=U$ and $M_\mathcal{S}=M_\text{tot}$.  $r_\mathcal{S}$ is the size of cluster $\mathcal{S}$.

We define the following functions:
\begin{align}
J_\mathcal{G}(q):=\sum_{i \in \mathcal{G}} m_i (q_i-L_i)^2 \quad \text{ and }\quad J_g(q):=\sum_{i=1}^n m_i (q_i-L_i)^2,  \label{eq:J_G-def}
\end{align}

The following theorem is the main result in this section.  It generalizes to a partial collision well known results for total collisions, see for example \cite{Wintner41}. Similar results  for partial collision has been established by Sperling \cite{Sp}, who  followed Wintner's approach (just as we do). Other result in this
direction is by Elbialy \cite{E}, who derived it using the McGehee regularization and exploiting the Shub lemma about compactness of the set of central configurations on the mass ellipsoid \cite{Sh70}.  We decided to include our  proof of this result for the sake of completeness.
\begin{theorem}
\label{thm:estm-coll-time}
Assume that bodies in cluster $\mathcal{G}$ collide for $t \nearrow T$.
Then, there exists a constant $A>0$ such that
\begin{align}
  \frac{J_{\mathcal{G}}(t)}{(T-t)^{4/3}} \to A, \quad t \nearrow T  \label{eq:thmJ-estm} \\
  \frac{\dot{J}_{\mathcal{G}}(t)}{(T-t)^{1/3}} \to -\frac{4A}{3}, \quad t \nearrow T \label{eq:thmdotJ-estm} \\
  \frac{\ddot{J}_\mathcal{G}(t)}{(T-t)^{-2/3}} \to \frac{4}{9}  A,  \quad t \nearrow T  \label{eq:ddotJ-bnds} \\
  \frac{U_\mathcal{G}(t)}{(T-t)^{-2/3}} \to \frac{2}{9}  A,  \quad t \nearrow T \label{eq:U-bnds} \\
  \frac{K_\mathcal{G}(t)}{(T-t)^{-2/3}} \to \frac{2}{9}  A,  \quad t \nearrow T. \label{eq:KG-bnds}
\end{align}
Moreover, the energy $H_\mathcal{G}(t)=K_\mathcal{G}(t) - U_\mathcal{G}(t)$ has a limit for $t \nearrow T $.
\end{theorem}

Before we proceed with the proof of Theorem~\ref{thm:estm-coll-time} we present two results about the behavior of the center of mass and the angular momentum
of cluster of bodies undergoing collision.

\begin{theorem}
\label{thm:cmass-beha}
Assume that $L_\mathcal{G}=0$. Then
for the center of mass $c_\mathcal{G}$, as $t$ approaches $T$, the following statements hold:
\begin{eqnarray}
  c_\mathcal{G}&=&O(|T-t|),  \label{eq:cG-asympt} \\
  \dot{c}_\mathcal{G}&=&O(1), \label{eq:dotcG-asympt} \\
   \ddot{c}_\mathcal{G}&=&O(1), \label{eq:ddotcG-asympt} \\
  M_\mathcal{G} c_\mathcal{G} \times \dot{c}_\mathcal{G}&=&O(|T-t|^2).  \label{eq:ang-mom-cmass}
\end{eqnarray}
\end{theorem}
\begin{proof}
Since
\begin{eqnarray*}
  M_\mathcal{G} \ddot{c}_\mathcal{G}=\sum_{i \in \mathcal{G}} \frac{\partial U_{\text{ext}}}{\partial q_i}
\end{eqnarray*}
and our trajectory is away of singularities in $U_{\text{ext}}$, it follows that $\ddot{c}_\mathcal{G}(t)$ is bounded for $t \in [t_0,T)$. Consequently, $\dot{c}_\mathcal{G}(t)$ converges as $t \nearrow T$, and since $c_\mathcal{G}(t) \to 0$ for $t \nearrow T$, we deduce (\ref{eq:cG-asympt}).

Furthermore, observe that
\begin{equation*}
  M_\mathcal{G} c_\mathcal{G}(t) \times \dot{c}_\mathcal{G}(t) \to 0, \quad t \nearrow T.
\end{equation*}
Hence, we obtain
\begin{eqnarray*}
  \frac{d}{dt} \left(M_\mathcal{G} c_\mathcal{G} \times \dot{c}_\mathcal{G}\right)=M_\mathcal{G} c_\mathcal{G} \times \ddot{c}_\mathcal{G} =O(|T-t|),
\end{eqnarray*}
and integrating yields (\ref{eq:ang-mom-cmass}).
\end{proof}

We reformulate estimates Theorem~\ref{thm:estm-coll-time} in terms of $r_\mathcal{G}^2=I^0_\mathcal{G}$.
\begin{theorem}
\label{thm:r-estm-coll-time}
Assume that bodies in cluster $\mathcal{G}$ collide for $t \nearrow T$.
Then, there exists a constant $A>0$ such that
\begin{align}
  \frac{I^0_{\mathcal{G}}(t)}{(T-t)^{4/3}} \to A, \quad t \nearrow T  \label{eq:r-thmJ-estm} \\
  \frac{\dot{I}^0_{\mathcal{G}}(t)}{(T-t)^{1/3}} \to -\frac{4A}{3}, \quad t \nearrow T \label{eq:r-thmdotJ-estm} \\
  \frac{\ddot{I}^0_\mathcal{G}(t)}{(T-t)^{-2/3}} \to \frac{4}{9}  A,  \quad t \nearrow T  \label{eq:r-ddotJ-bnds} \\
\end{align}
\end{theorem}
\begin{proof}
Without any loss of generality we can assume that $L_\mathcal{G}=0$. Then $J_\mathcal{G}(t)=I^0_\mathcal{G}(t) + M_\mathcal{G} c^2_\mathcal{G}(t)$ and  from (\ref{eq:thmJ-estm},\ref{eq:thmdotJ-estm},\ref{eq:ddotJ-bnds}) in Theorem~\ref{thm:estm-coll-time} and Theorem~\ref{thm:cmass-beha}
we obtain
\begin{eqnarray*}
  \frac{I^0_{\mathcal{G}}(t)}{(T-t)^{4/3}} =  \frac{J_{\mathcal{G}}(t) + O(|T-t|^2)}{(T-t)^{4/3}}  \to A, \quad t \nearrow T   \\
  \frac{\dot{I}^0_{\mathcal{G}}(t)}{(T-t)^{1/3}} = \frac{\dot{J}_{\mathcal{G}}(t) + O(|T-t|)}{(T-t)^{1/3}} \to -\frac{4A}{3}, \quad t \nearrow T  \\
   \frac{\ddot{I}^0_\mathcal{G}(t)}{(T-t)^{-2/3}}=\frac{\ddot{J}_\mathcal{G}(t) + O(1)}{(T-t)^{-2/3}} \to \frac{4}{9}  A,  \quad t \nearrow T.
\end{eqnarray*}
\end{proof}

\begin{lemma}
\label{lem:r-ver-JG}
Assume that bodies in cluster $\mathcal{G}$ collide for $t \nearrow T$.
 Then
\begin{equation}
  \frac{r_\mathcal{G}(t)}{\sqrt{J_\mathcal{G}(t)}} \to 1, \quad t \nearrow T.
\end{equation}
\end{lemma}
\begin{proof}
Let us denote $q_\mathcal{G}=(q_i)_{i \in \mathcal{G}}$.
From the triangle inequality applied to the norm $\|x\|=\sqrt{\sum_{i\in \mathcal{G}}m_i x_i^2}$  it follows that
\begin{eqnarray*}
   \|q_\mathcal{G}-L_\mathcal{G}\| -  \|L_\mathcal{G}-c_\mathcal{G}\| \leq  \|q_\mathcal{G}-c_\mathcal{G}\| \leq  \|q_\mathcal{G}-L_\mathcal{G}\| +  \|L_\mathcal{G}-c_\mathcal{G}\|,
\end{eqnarray*}
which written in terms of $r_\mathcal{G}$ and $J_\mathcal{G}$ becomes
\begin{equation*}
  \sqrt{J_\mathcal{G}} - \|L_\mathcal{G}-c_\mathcal{G}\| \leq r_\mathcal{G} \leq  \sqrt{J_\mathcal{G}} + \|L_\mathcal{G}-c_\mathcal{G}\|,
\end{equation*}
and after division by $\sqrt{J_\mathcal{G}}$ we obtain
\begin{equation*}
  1- \frac{\|L_\mathcal{G}-c_\mathcal{G}\|}{\sqrt{J_\mathcal{G}}} \leq \frac{r_\mathcal{G}}{\sqrt{J_\mathcal{G}}} \leq   1+ \frac{\|L_\mathcal{G}-c_\mathcal{G}\|}{\sqrt{J_\mathcal{G}}},
\end{equation*}
Since from Theorems~\ref{thm:estm-coll-time} and \ref{thm:cmass-beha} $ \frac{\|L_\mathcal{G}-c_\mathcal{G}\|}{\sqrt{J_\mathcal{G}}}=O(|T-t|^{1/3})$ we obtain our assertion.
\end{proof}

The next theorem due to Elbialy \cite{E} gives bounds for the intrinsic angular momentum of the cluster.   We present here its proof following ideas of Elbialy for the sake
of completeness.
\begin{theorem}\cite[Proposition 4.11]{E}
\label{thm:intr-angmom-collison}
Assume that bodies in cluster $\mathcal{G}$ collide for $t \nearrow T$. Then,
the intrinsic angular momentum of bodies in cluster $\mathcal{G}$ satisfies
\begin{eqnarray}
  \mu^0_{\mathcal{G}}(t) &=& O(|T-t|^{\frac{7}{3}}) \label{eq:muint}  \\
  \dot{\mu}^0_{\mathcal{G}}(t)&=&O(|T-t|^{\frac{4}{3}}),
\end{eqnarray}
\end{theorem}
\begin{proof}
Let us denote $\mu^0_{\mathcal{G}}$ by $\mu$.  We have (compare (\ref{eq:mu-corr-cart}))
\begin{equation*}
 \mu= \sum_{i\in\mathcal{G}} m_i q_i \times \dot{q}_i -  M_{\mathcal{G}} c_{\mathcal{G}}  \times \dot{c}_{\mathcal{G}}=\mu_\mathcal{G} - M_{\mathcal{G}} c_{\mathcal{G}}\times \dot{c}_{\mathcal{G}}.
\end{equation*}
Without any loss of generality we can assume that $L_\mathcal{G}=0$. From  (\ref{eq:ang-mom-cmass}) in Theorem~\ref{thm:cmass-beha}
we have
\begin{equation}
 M_{\mathcal{G}} c_{\mathcal{G}}\times \dot{c}_{\mathcal{G}}=O(|T-t|^2). \label{eq:amomcm}
\end{equation}

From the asymptotic bounds stated in Theorem~\ref{thm:estm-coll-time}, we have
\begin{equation*}
  K_\mathcal{G}(t) = O(|T-t|^{-2/3}),
\end{equation*}
which implies bounds on velocities   $|\dot{q}_i|=O(|T-t|^{-1/3})$ for $i \in \mathcal{G}$. Consequently,
\begin{eqnarray*}
 |\mu_\mathcal{G}(t)| \leq O(|T-t|^{2/3}) \cdot O(|T-t|^{-1/3})= O(|T-t|^{1/3}),
\end{eqnarray*}
from this and (\ref{eq:amomcm}) we obtain
\begin{equation}
  \lim_{t \nearrow T} \mu(t)=0.  \label{eq:mu(T)=0}
\end{equation}

Now we estimate the derivative of $\mu$. It is easy to see that
\begin{eqnarray*}
  \dot{\mu}=\sum_{i \in \mathcal{G}} m_i (q_i - c_\mathcal{G}) \times (\ddot{q}_i - \ddot{c}_\mathcal{G}) = \sum_{i \in \mathcal{G}} m_i (q_i - c_\mathcal{G}) \times \ddot{q}_i  \\
  =  \sum_{i \in \mathcal{G}} (q_i - c_\mathcal{G}) \times \left( \sum_{j \notin \mathcal{G}} m_i m_j \frac{q_j-q_i}{|q_j - q_i|^3}\right)
\end{eqnarray*}
Each of the terms $ \frac{q_j-q_i}{|q_j - q_i|^3}$ we can write as follows (this Taylor formula with 1-st order remainder)
\begin{eqnarray*}
   \frac{q_j-q_i}{|q_j - q_i|^3}= \frac{q_j-c_\mathcal{G}}{|q_j - c_\mathcal{G}|^3} + b(q_j,q_i,c_\mathcal{G})(q_i - c_\mathcal{G})
\end{eqnarray*}
where $b()$ is some bounded matrix.

Observe that for each $j \notin \mathcal{G}$ holds $ \sum_{i \in \mathcal{G}} (q_i - c_\mathcal{G}) \times (\frac{q_j-c_\mathcal{G}}{|q_j - c_\mathcal{G}|^3})=0$, while
the first order remainder term times $q_i - c_\mathcal{G}$  is quadratic in $|q_i-c_\mathcal{G}|$, hence is $O(r^2)$. Therefore we obtained
\begin{equation}
  \dot{\mu}(t)=O(r^2(t))).
\end{equation}
Now from Theorem~\ref{thm:estm-coll-time} and Lemma~\ref{lem:r-ver-JG} we obtain
\begin{equation}
\dot{\mu}(t)=O(|T-t|^{4/3}).
\end{equation}
After integration of the above equation and  taking into account (\ref{eq:mu(T)=0}) we obtain (\ref{eq:muint}).
\end{proof}

\subsection{Proof of Theorem~\ref{thm:estm-coll-time}}

To establish Theorem~\ref{thm:estm-coll-time}, we adapt the arguments from Wintner's book \cite{Wintner41}, where the total collision was considered.

In Pollard and Saari paper\cite{PS68}, where the partial collision was considered, the bounds are obtained for $J_g$, $\dot{J}_g$ and $U$ and not separately for clusters of colliding bodies, i.e., $J_\mathcal{G}$, $\dot{J}_\mathcal{G}$ and $U_\mathcal{G}$.
Furthermore, the bounds obtained by \cite{PS68} rely on a particular Tauberian theorem, which is essentially based on a proposition attributed to Boas \cite[Thm 1]{B39}. Upon closer examination, we have reservations about some steps in the proof of this proposition.

The proof of Theorem \ref{thm:estm-coll-time}  is organized as follows:  first, we establish (\ref{eq:thmJ-estm}) and (\ref{eq:thmdotJ-estm}) for $J_g$ - this the content of Lemma~\ref{lem:estm-collJg}. From this we infer that $H_\mathcal{G}$ is bounded (see Lemma~\ref{lem:HC-bnd}),
which allows us to repeat the same reasoning which led to the result for $J_g$ also to $J_\mathcal{G}$ and to obtain  (\ref{eq:thmJ-estm}) and (\ref{eq:thmdotJ-estm}). This is
the content of Subsection~\ref{sssec:firstpart}. Then in Subsection~\ref{sssec:sec-part} we establish remaining conditions (\ref{eq:ddotJ-bnds},\ref{eq:U-bnds},\ref{eq:KG-bnds}).

\subsubsection{First part - proof of  (\ref{eq:thmJ-estm}) and (\ref{eq:thmdotJ-estm})}
\label{sssec:firstpart}
\begin{lemma}
\label{lem:estm-collJg}
Under assumptions of Theorem~\ref{thm:estm-coll-time} there exists $A>0$ such that
\begin{eqnarray}
 \frac{J_{g}(t)}{(T-t)^{4/3}} \to A, \quad t \nearrow T  \label{eq:thmJg-estm} \\
 \frac{\dot{J}_g(t)}{(T-t)^{1/3}} \to -\frac{4A}{3}, \quad t \nearrow T.  \label{eq:thmdotJg-estm}
\end{eqnarray}
\end{lemma}

To obtain the assertion regarding $J_\mathcal{G}$, the argument is almost the same as in the proof of Lemma~\ref{lem:estm-collJg}, but we use additional bounds obtained from (\ref{eq:thmdotJg-estm}).
Consequently, we derive several lemmas that apply to both $J_g$ and $J_{\mathcal{G}}$.
However, an additional assumption arises, which we can verify once we obtain the bounds from Lemma \ref{lem:estm-collJg} on $J_g$. Subsequently, we address (\ref{eq:U-bnds}).

We begin by assuming that the center of mass is located at the origin, that is,
\begin{eqnarray}
\sum_i m_i q_i=0,   \label{eq:ccm=0}
\end{eqnarray}
 this also implies that
\begin{eqnarray}
\sum_i m_i L_i=0.  \label{eq:ccL}
\end{eqnarray}
Observe that $J_{\mathcal{G}}$ and $J_g$ do not depend upon assumption (\ref{eq:ccm=0}) (or (\ref{eq:ccL})).
These are technical assumptions that help and, on the other side, have a reasonable physical meaning.

In what follows, we will apply the same reasoning for $J_g$ and $J_{\mathcal{G}}$, denoting them simply as $J$.
When $J=J_g$, we will utilize $\mathcal{P}=\{1,\dots,n\}$, and when considering $J=J_{\mathcal{G}}$, we will use $\mathcal{P}=\mathcal{G}$.

From Wintner's proof, it is clear that the goal is to prove that
\begin{equation*}
  \frac{\dot{J}^2}{\sqrt{J}} \to b >0, \quad \dot{J}<0, \qquad t \nearrow T.
\end{equation*}
This is a content of Lemma~\ref{lem:limQ}.
From this, it is relatively straightforward to derive bounds (\ref{eq:thmJ-estm},\ref{eq:thmdotJ-estm}).

Since $U$ is homogeneous of degree $-1$, we have that
\begin{equation}
  \sum_{i \in \mathcal{P}} (q_i, \nabla_{q_i} U_\mathcal{P})=-U_\mathcal{P},  \label{eq:qgradU}
\end{equation}
and due to the translational invariance of $U_\mathcal{P}$, we have
\begin{equation}
   \sum_{i \in \mathcal{P}} \nabla_{q_i} U_\mathcal{P}= 0. \label{eq:sumgrad=0}
\end{equation}

\begin{lemma}
\label{lem:ddj}
\begin{eqnarray}
  \ddot{J}_\mathcal{P}(t)&=& 4 K_\mathcal{P}(t) - 2U_\mathcal{P}(t) + g(t) = 4H_\mathcal{P}(t) + 2 U_\mathcal{P}(t) + g(t) =2 K_{\mathcal{P}}(t) + 2H_\mathcal{P}(t) + g_\mathcal{P}(t),  \nonumber
\end{eqnarray}
where $g_\mathcal{P}(t)$ is a bounded continuous function.
\end{lemma}
\begin{proof}
 We have
\begin{eqnarray*}
  \dot{J}&=&\sum_{i \in \mathcal{P}} 2 m_i (\dot{q}_i,q_i-L_i), \\
  \ddot{J}&=& \sum_{i \in \mathcal{P}} 2 m_i (\dot{q}_i,\dot{q}_i) + \sum_{i\in \mathcal{P}} 2 m_i  (\ddot{q}_i,q_i-L_i)= 4K_\mathcal{P} + 2 \sum_{i\in \mathcal{P}} (\nabla_i U(q),q_i - L_i).
\end{eqnarray*}
Let us take a look at $\sum_{i \in \mathcal{P}} (\nabla_i U(q),q_i - L_i)$.
We split it into intra-cluster and inter-cluster interactions.
It will turn out that the intra-cluster contribution, which contains unbounded terms, vanishes and what remains is a bounded function of $t \in [t_0,T)$.

We have
\begin{equation*}
  U(q) = U_{\text{in}}(q) + U_{\text{ex}}(q),
\end{equation*}
where
\begin{equation*}
 U_{\text{in}}(q) =\sum_{k=1,\dots,N_c} U_{G_k}, \quad\quad
 U_{\text{ex}}(q) =\sum_{j<k, \mbox{not in the same cluster}} U_{k,j}(q_k,q_j).
\end{equation*}
Observe that $U_{\text{ex}}(q)$ is a very regular function in the region where the collision occurs.

Now, if $J=J_{\mathcal{G}}$, then  by applying (\ref{eq:qgradU},\ref{eq:sumgrad=0}), we have
\begin{eqnarray*}
  \sum_{i \in \mathcal{G}} (\nabla_i U(q),q_i - L_i) = \sum_{i \in \mathcal{G}}\left( (\nabla_i U_\mathcal{G}(q),q_i) - (\nabla_i U_\mathcal{G}(q),L_\mathcal{G})
     + (\nabla_i U_{\text{ex}}(q),q_i-L_\mathcal{G}))\right) \\
   = -U_\mathcal{G}(q) - \left(\sum_{i \in \mathcal{G}}\nabla_i U_\mathcal{G}(q) ,L_\mathcal{G}\right) + g(t)  = -U_\mathcal{G}(q) + g(t)
\end{eqnarray*}
where $g(t)=(\nabla_i U_{\text{ex}}(q),q_i-L_\mathcal{G}))$ is continuous bounded function for $t \in [t_0,T)$.
Therefore, we obtain
\begin{eqnarray*}
  \ddot{J}= 4K_\mathcal{P} - 2U_\mathcal{P} + 2 g(t)=4 H_\mathcal{P} + 2U_\mathcal{P} + 2 g(t).
\end{eqnarray*}
This is our assertion.

Now, we consider the case of $J=J_g$, i.e., $\mathcal{P}=\{1,\dots,n\}$. From (\ref{eq:qgradU}), we have
\begin{eqnarray*}
\sum_{i \in \mathcal{P}} (\nabla_i U(q),q_i - L_i) = \sum_{i} (\nabla_i U(q),q_i) - \sum_i (\nabla_i U(q),L_i) = -U(q) - \sum_i (\nabla_i U(q),L_i).
\end{eqnarray*}

We have from (\ref{eq:sumgrad=0})
\begin{eqnarray*}
 \sum_i (\nabla_i U(q),L_i) =  \sum_{k} \sum_{i \in G_k} (\nabla_i U(q),L_{G_k})\\
  =  \sum_{k} \sum_{i \in G_k} (\nabla_i U_{G_k}(q),L_{G_k})  + \sum_{k} \sum_{i \in G_k} (\nabla_i U_{\text{ex}}(q),L_{G_k}) \\
  =  \sum_{k} \left(\sum_{i \in G_k}\nabla_i U_{G_k}(q),L_{G_k}\right) + g(t)=g(t)
\end{eqnarray*}
where $g(t)=\sum_{i \in \mathcal{P}} (\nabla_i U_{\text{ex}}(q),L_i)$.

We proceed as in the previous case.

\end{proof}

From Lemma~\ref{lem:ddj}, we obtain the following result.
\begin{lemma}
\label{lem:dJ<0}
Under the same assumptions as in Thm~\ref{thm:estm-coll-time}.
In the case of $J=J_\mathcal{G}$, assuming additionally that $H_\mathcal{G}(t)$ is bounded for $t \in [t_0,T)$, $\dot{J}<0$  for $t$ close enough to $T$.
\end{lemma}

\begin{proof}
From Lemma~\ref{lem:ddj}
\begin{equation*}
  \ddot{J}_\mathcal{P}=4H_\mathcal{P}(t) + 2 U_\mathcal{P}(t) + g(t)= 2U_\mathcal{G}(t) + (4H_\mathcal{G}(t) + g(t))=2U_\mathcal{G}(t) + \tilde{g}(t).
\end{equation*}
where $\tilde{g}(t)$ is a bounded function, because if $J=J_g$, then $H_\mathcal{P}(t)=h$ by the conservation of energy and for $J=J_\mathcal{G}$ we assumed that $H_\mathcal{G}(t)$ is bounded.

Therefore, we have
\begin{equation*}
  \ddot{J}(t)=2 U_\mathcal{P}(t) + \tilde{g(t)}.
\end{equation*}
Since $U_\mathcal{P}(t) \to \infty$ as $t\nearrow T$, $\ddot{J} >0$ for $t \geq t_0$, hence $\dot{J}$ is an increasing function for $t>t_0$. Therefore, there exists $t_1 \geq t_0$ such that $\dot{J} <0$, as otherwise
$J(t)>0$ would be non-decreasing, contradicting the convergence $J(t) \to 0$ as $t \nearrow T$.
\end{proof}

Next three lemma are inequalities involving $J_\mathcal{P}$, $K_\mathcal{P}$ and $U_\mathcal{P}$.
\begin{lemma}
\label{lem:dJ<2jk}
For any subset $\mathcal{P} \subset \{1,\dots,n\}$, the following inequality holds
 \begin{equation*}
    \frac{\dot{J}_\mathcal{P}^2}{4} \leq 2 J_\mathcal{P} K_\mathcal{P}.
 \end{equation*}
\end{lemma}

\begin{proof}
We apply the Schwarz inequality to the scalar product $\langle x,y \rangle=\sum_i m_i (x_i,y_i)$.
Thus,
\begin{eqnarray*}
  \left(\frac{\dot{J}_\mathcal{P}}{2}\right)^2= \left( \sum_{i \in \mathcal{P}} m_i (\dot{q}_i,q_i-L_i) \right)^2 \leq \left(\sum_{i \in \mathcal{P}} m_i \dot{q}_i^2 \right)  \left(\sum_{i \in \mathcal{P}} m_i (q_i - L_i)^2 \right) \leq 2J_\mathcal{P}K_\mathcal{P}.
\end{eqnarray*}
\end{proof}

\begin{lemma}
\label{lem:J-mutual-dist}
The following equalities are satisfied
  \begin{equation*}
    J_g=M_\text{\rm tot}^{-1} \sum_{j < k} m_j m_k \|(q_j - L_j) - (q_k - L_k) \|^2
  \end{equation*}
 and
 \begin{equation*}
   J_{\mathcal{G}}=M_\mathcal{G}^{-1} \sum_{j < k, j,k \in \mathcal{G}} m_j m_k \|(q_j - L_j) - (q_k - L_k) \|^2 + M_\mathcal{G}\|c_\mathcal{G} - L_\mathcal{G}\|^2.
 \end{equation*}
\end{lemma}

\begin{proof}
Let us set $\hat{q}_k = q_k - L_k$.
Then,
\begin{eqnarray*}
  2 \sum_{j < k, j,k \in \mathcal{P}} m_j m_k \| \hat{q}_j - \hat{q}_k \|^2 = \sum_{ j,k \in \mathcal{P}} m_j m_k \| \hat{q}_j - \hat{q}_k \|^2 =
  \sum_{ j,k \in \mathcal{P}} m_j m_k  \left(\hat{q}_j^2 -2 (\hat{q}_k,\hat{q}_j) +  \hat{q}_k^2 \right) \\
  = 2 M_\mathcal{P} J_\mathcal{P} - 2 \sum_{k \in \mathcal{P}}  m_k \left(\hat{q}_k, \sum_{j \in \mathcal{P}} m_j \hat{q}_j \right).
\end{eqnarray*}

If $\mathcal{P}=\{1,\dots,n\}$, then from \eqref{eq:ccm=0} and \eqref{eq:ccL}, we have
\begin{eqnarray*}
 \sum_{j \in \mathcal{P}} m_j \hat{q}_j=  \sum_{j} m_j q_j   - \sum_{j} m_j L_j =0.
\end{eqnarray*}
Therefore, the assertion is proved.

If $\mathcal{P}=\mathcal{G}$, then $L_j=L_\mathcal{G}$ and we obtain
\begin{eqnarray*}
  \sum_{ j,k \in \mathcal{G}} m_j m_k  (\hat{q}_k,\hat{q}_j) = \left( \sum_{k \in \mathcal{G}} m_k (q_k - L_\mathcal{G}) , \sum_{j \in \mathcal{G}} m_j (q_j - L_\mathcal{G})\right) =
   \left( M_\mathcal{G} (c_\mathcal{G} - L_\mathcal{G}) ,  M_\mathcal{G} (c_\mathcal{G} - L_{\mathcal{G}})\right).
\end{eqnarray*}
Thus, we obtain
\begin{eqnarray*}
   2 \sum_{j < k, j,k \in \mathcal{G}} m_j m_k \| \hat{q}_j - \hat{q}_k \|^2 =2 M_\mathcal{G} J_{\mathcal{G}}  - 2 M_\mathcal{G}^2 \|c_\mathcal{G} - L_\mathcal{G}\|^2.
\end{eqnarray*}
\end{proof}

\begin{lemma}
\label{lem:sqrtJU>D}
There exists a positive constant $D$, depending on the masses, such that
\begin{equation*}
  \sqrt{J_\mathcal{P}} U_\mathcal{P} \geq D.
\end{equation*}
\end{lemma}

\begin{proof}
Let's define
\begin{equation*}
  J_{\text{rel}}=M_\mathcal{P}^{-1} \sum_{j < k, j,k \in \mathcal{P}} m_j m_k \|(q_j - L_j) - (q_k - L_k)\|^2.
\end{equation*}
From Lemma~\ref{lem:J-mutual-dist} it follows that $J_\mathcal{P} \geq J_{\text{rel}}$.

For indices $j,k \in \mathcal{P}$ in the same colliding cluster, i.e., $L_j=L_k$, it holds that $ \|(q_j - L_j) - (q_k - L_k)\|= \|q_j - q_k\|$.
Therefore,
\begin{eqnarray*}
  \sqrt{J_{\text{rel}}} \geq \sqrt{M_\mathcal{P}^{-1} m_j m_k} \cdot \|q_{j}-q_k\|,
\end{eqnarray*}
hence
\begin{equation*}
   \sqrt{J_{\text{rel}}} \frac{m_j m_k}{\|q_{j}-q_k\|} \geq M_\mathcal{P}^{-1/2} m_j^{3/2} m_k^{3/2}.
\end{equation*}
Finally,
\begin{eqnarray*}
  \sqrt{J_{\text{rel}}} U_{\mathcal{P}}  \geq  M_\mathcal{P}^{-1/2} \sum_{j < k, j,k \in \mathcal{P}}m_j^{3/2} m_k^{3/2}=D.
\end{eqnarray*}
Given that  $J_{\mathcal{P}} \geq J_{\text{rel}}$ the assertion follows.
\end{proof}

Now, we are ready for the crucial lemma.
\begin{lemma}
\label{lem:limQ}
Under the same assumptions as in Thm~\ref{thm:estm-coll-time}. In the case of $J=J_\mathcal{G}$, we assume additionally that $H_\mathcal{G}(t)$ is bounded for $t \in [t_0,T)$.

 There exists $b>0$, such that
  \begin{equation*}
    \lim_{t \nearrow T}\frac{\dot{J}^2}{J^{1/2}} = b.
  \end{equation*}
\end{lemma}

\begin{proof}
    Let us set
\begin{equation*}
  Q(t)=A J^{1/2} + \frac{\dot{J}^2}{J^{1/2}}, \label{eq:def-Q}
\end{equation*}
where $A$ is a suitable constant to be determined later.
We aim to prove that $\dot{Q}(t) <0$ near the collision.

In the reasoning below we work with $t \in [t_0,T)$ such that $\dot{J}(t) <0$ (see Lemma~\ref{lem:dJ<0}). We have a freedom to increase $t_0$ if needed.

We have
\begin{eqnarray*}
\dot{Q}=  A \dot{\sqrt{J}} + 4 \dot{\sqrt{J}}\left(\ddot{J} - \frac{\dot{J}^2}{4J} \right)= 4 \dot{\sqrt{J}} \left(\frac{A}{4} + \ddot{J} - \frac{\dot{J}^2}{4J} \right).
\end{eqnarray*}
We want $A$ such that $\left(\frac{A}{4} + \ddot{J} - \frac{\dot{J}^2}{4J} \right)>0$ for $t$ close to $T$.
From Lemmas~\ref{lem:ddj} and~\ref{lem:dJ<2jk}, we have
\begin{eqnarray*}
  \frac{A}{4} + \ddot{J} - \frac{\dot{J}^2}{4J} \geq \frac{A}{4} +( 4K_\mathcal{P}(t) - 2U_\mathcal{P}(t) + g(t)) - 2K_\mathcal{P}(t) = 2H_\mathcal{P}(t) + g(t) + \frac{A}{4}.
\end{eqnarray*}
Therefore, since $\dot{\sqrt{J}}<0$ choosing
\begin{equation*}
  A > -4g(t) - 8H_{\mathcal{P}}(t),
\end{equation*}
ensures $\dot{Q} <0$.
Let us fix any such $A$; this is possible because $g(t)$ and $H_\mathcal{P}(t)$ are bounded functions; recall that $H_\mathcal{P}(t)=h$ if $J=J_g$.

$Q$ is monotonically decreasing and bounded from below by $0$, implying that $\lim_{t \nearrow T} Q(t)=b \geq 0$ exists.
Therefore, $b=\lim_{t\nearrow T} \frac{\dot{J}^2}{J^{1/2}}$.
Consequently, it follows that
\begin{equation*}
  \lim_{t \nearrow T} \dot{J}(t)^2=0.
\end{equation*}

Our goal is to show that $b>0$.
Since
\begin{eqnarray*}
  \frac{d}{dt} \dot{J}^2= 2 \dot{J} \ddot{J},
\end{eqnarray*}
by integrating both sides of the above equality, we obtain
\begin{eqnarray*}
 - \dot{J}(t)^2 = \dot{J}(T)^2  - \dot{J}(t)^2   = 2\int_t^T \dot{J}(s) \ddot{J}(s)ds.
\end{eqnarray*}
From Lemmas~\ref{lem:ddj} and~\ref{lem:sqrtJU>D}, we have for some $D>0$
\begin{eqnarray*}
  \ddot{J} =  4H_{\mathcal{P}}(t) + 2U_\mathcal{P}(t) + g(t)  \geq \left( 4 H_\mathcal{P}(t) + g(t)\right) + \frac{2D}{\sqrt{J}}.
\end{eqnarray*}

Since $J(t) \to 0$, for $t$ sufficiently close to $T$, we will have
\begin{equation}
    \ddot{J} > \frac{D}{ \sqrt{J}}.  \label{eq:ex-under-int}
\end{equation}

Considering the fact that  $\dot{J}<0$  and from \eqref{eq:ex-under-int}, for $t$ close enough to $T$,
\begin{eqnarray*}
  -\int_t^T \dot{J}(s) \ddot{J}(s)ds= \int_t^T (-\dot{J}(s)) \ddot{J}(s)ds >  \int_t^T (-\dot{J}(s))\frac{D}{\sqrt{J(s)}}ds\\
  = -2 D \int_t^T \dot{\sqrt{J}} ds = 2D \sqrt{J(t)} - 2D \sqrt{J(T)}=  2D \sqrt{J(t)}.
\end{eqnarray*}
Combining the above computation, we obtain for $t$ close enough to $T$,
\begin{eqnarray*}
  \dot{J}(t)^2 = - 2\int_t^T \dot{J}(s)\ddot{J}(s)ds \geq 4D \sqrt{J(t)}.
\end{eqnarray*}
Hence
\begin{equation*}
 \frac{\dot{J}^2}{\sqrt{J}} \geq 4D.
\end{equation*}
Therefore
\begin{equation*}
  b \geq 4D >0.
\end{equation*}
\end{proof}

\begin{lemma}
\label{lem:JdotJ-estm}
Under the same assumptions as in Thm~\ref{thm:estm-coll-time}. In the case of $J=J_\mathcal{G}$, we assume additionally that $H_\mathcal{G}(t)$ is bounded for $t \in [t_0,T)$.

Then
  \begin{eqnarray}
     \frac{J(t)}{(T - t)^{4/3}} \to A, \quad t \nearrow T , \nonumber\\
     \frac{\dot{J}(t)}{(T - t)^{1/3}} \to   -\frac{4A}{3}, \quad t \nearrow T, \nonumber
  \end{eqnarray}
  where $A=\left(\frac{3\sqrt{b}}{4} \right)^{4/3}$ and $b$ is a constant obtained in Lemma~\ref{lem:limQ}.
\end{lemma}

\begin{proof}
From Lemma~\ref{lem:dJ<0} we can assume that $\dot{J}<0$.

    From Lemma~\ref{lem:limQ}, it follows that for any $\varepsilon>0$ for $t$ close enough to $T$, the following inequality holds
\begin{equation*}
   \left| \frac{-\dot{J}(t)}{\sqrt{b} J(t)^{1/4}} - 1  \right| < \varepsilon.
\end{equation*}
Given that $\dot{J} <0$, we can rewrite this as
\begin{equation}
 -(1+\varepsilon) \sqrt{b}  < \frac{\dot{J}}{ J^{1/4} } < - (1-\varepsilon) \sqrt{b}. \label{eq:dotJ-estm}
\end{equation}
Now, by defining the inverse function $t(J)$, which exists because $\dot{J}<0$ for $t$ close to $T$, we have
\begin{align*}
      -(1-\varepsilon)\sqrt{b}\frac{dt}{dJ} &<  \frac{1}{J^{1/4}} < -(1+\varepsilon)\sqrt{b_0}\frac{dt}{dJ}.
\end{align*}
Integrating the above inequality from 0 to $J$ yields:
\begin{align*}
     (1-\varepsilon)\sqrt{b}\int_0^J\left(-\frac{dt}{dJ}\right)dJ &<\int_0^J  \frac{dJ}{J^{1/4}} < (1+\varepsilon)\sqrt{b}\int_0^J\left(-\frac{dt}{dJ}\right)dJ\\
    (1-\varepsilon)\sqrt{b}(T-t) &< \frac{4}{3}J^{3/4} < (1+\varepsilon)\sqrt{b}(T-t)\\
    \left[\frac{3(1-\varepsilon)\sqrt{b}}{4}\right]^{4/3}(T - t)^{4/3} & < J < \left[\frac{3(1+\varepsilon)\sqrt{b}}{4} \right]^{4/3}(T - t)^{4/3}.\\
\end{align*}
Hence, $A=\left(\frac{3\sqrt{b}}{4} \right)^{4/3}$.

From the above inequality and \eqref{eq:dotJ-estm}, we obtain
\begin{eqnarray*}
 -(1+\varepsilon) \sqrt{b} \left[\frac{3(1+\varepsilon)\sqrt{b}}{4} \right]^{1/3}(T - t)^{1/3}  \leq    \dot{J} \leq - (1-\varepsilon) \sqrt{b}  \left[\frac{3(1-\varepsilon)\sqrt{b}}{4}\right]^{1/3}(T - t)^{1/3}.
\end{eqnarray*}
\end{proof}

To finish the proof  (\ref{eq:thmJ-estm},\ref{eq:thmdotJ-estm}) in Theorem~\ref{thm:estm-coll-time}, we must show the following lemma.

\begin{lemma}
\label{lem:HC-bnd}
Under the same assumptions as in Thm~\ref{thm:estm-coll-time}, the function
\begin{equation}
  H_\mathcal{G}(t)=K_\mathcal{G}(t) - U_\mathcal{G}(t), \label{eq:Hc-bounded}
\end{equation}
has a limit $\lim_{t \nearrow T} H_\mathcal{G}(t)$. In particular,  $H_\mathcal{G}(t)$ is bounded for $t \in [t_0,T)$.
\end{lemma}

\begin{proof}
Recall, from \eqref{eq:second-order-formulation} and \eqref{eq:energy_potential_inCluster}, that the equation of motion of a particle in cluster $\mathcal{G}$ is
\begin{eqnarray*}
  m_i\ddot{q}_i = \frac{\partial U_\mathcal{G}}{\partial q_i} +  \frac{\partial U_{\text{ext}}}{\partial q_i}, \quad i \in \mathcal{G}.
\end{eqnarray*}
Then, the rate of change of internal energy for the cluster will be given by
\begin{eqnarray}
  \frac{d}{dt} H_\mathcal{G}(t) = \sum_{i \in \mathcal{G}} \frac{\partial U_{\text{ext}}}{\partial q_i} \cdot \dot{q}_i(t).  \label{eq:dE-per}
\end{eqnarray}
While the terms $\frac{\partial U_{\text{ext}}}{\partial q_i}$ are bounded, $\dot{q}_i(t)$ may not be.
However, since the integration time is finite,  we expect that the integrated change of $\int_{t_0}^T \frac{dH_\mathcal{G}}{dt}(t)dt$ will be bounded.
This holds true when $q(t)$ has a finite length.

We will use our estimates on $J_g$ to show that $q(t)$ has a finite length when approaching the collision, implying that the integral of $dH_\mathcal{G}/dt$ given by \eqref{eq:dE-per} is bounded.
For this end, we establish a bound on $\int_{t_0}^T K_\mathcal{G}(t)dt $, which, with the aid of the Schwarz inequality, leads to the desired result.

From Lemma~\ref{lem:JdotJ-estm}, for $J_g$, we obtain
\begin{equation}
  -  B_1\cdot (T - t)^{1/3}  < \dot{J}_g(t) <  -A_1\cdot (T - t)^{1/3}.  \label{eq:bnddotJg}
\end{equation}

From Lemma~\ref{lem:ddj}, we have
\begin{equation}
  \ddot{J}_g=4h + 2U +g(t), \label{eq:ddJg2}
\end{equation}
where
\begin{equation*}
 |g(t)| \leq B_g, \quad t  \in [t_0,T).
\end{equation*}

Thanks to Lemma~\ref{lem:dJ<0}, we can take $t_0$ such that $\ddot{J}_g(t)>0$ and $\dot{J}_g(t)<0$ for $t \in [t_0,T)$.
Integrating (\ref{eq:ddJg2}) from $t_0$ to $t \nearrow T$, we obtain for suitable constant $B_2$
\begin{align*}
  \dot{J}_g(t) - \dot{J}_g(t_0) &= 4(t-t_0)h + 2\int_{t_0}^t U(q(s)) ds + \int_{t_0}^t g(s) ds, \label{eq:eIntExpJ'}\\
  \int_{t_0}^t U(q(s)) ds &= \frac{1}{2}\left( \dot{J}_g(t) - \dot{J}_g(t_0) \right) -2 (t-t_0)h - \frac{1}{2}\int_{t_0}^t g(s) ds, \nonumber \\
  \int_{t_0}^t U(q(s)) ds  &\leq |\dot{J}_g(t_0)| + (T-t_0) \left(2|h| + B_g/2\right) \nonumber \\
      &\leq B_1 (T-t_0)^{2/3} +  (T-t_0) \left(2|h| + B_g/2\right)\leq B_2(T-t_0)^{2/3}. \nonumber
\end{align*}
Since $K=h+U$, we obtain, for some constant $B_3$,
\begin{equation*}
  \int_{t_0}^T \frac{1}{2}\sum_i m_i \dot{q}_i^2 = (T-t_0) h +  B_2(T-t_0)^{2/3} \leq B_3 (T-t_0)^{2/3}.
\end{equation*}

Returning to \eqref{eq:dE-per}, we find
\begin{eqnarray*}
  H_\mathcal{G}(t) - H_\mathcal{G}(t_0) =  \int_{t_0}^t \sum_{i \in \mathcal{G}} \frac{\partial U_{\text{ext}}}{\partial q_i} \cdot \dot{q}_i(s)ds.
\end{eqnarray*}
We know that there exists a constant $Z$ such that each term $\left|\frac{\partial U_{\text{ext}}}{\partial q_i} \right| \leq Z$.
Then, from Schwarz inequality, for some constant $\tilde{Z}$ for all $t \in [t_0,T)$, we have
\begin{eqnarray*}
  \left|\int_{t_0}^t \sum_{i \in \mathcal{G}} \frac{\partial U_{\text{ext}}}{\partial q_i} \cdot \dot{q}_i(s)ds \right|
  \leq \int_{t_0}^T \left| \sum_{i \in \mathcal{G}} \frac{\partial U_{\text{ext}}}{\partial q_i} \cdot \dot{q}_i(s)\right| ds  \leq
       \int_{t_0}^T \sum_{i \in \mathcal{G}}n Z |\dot{q}_i(s)| ds \\
        \leq \tilde{Z} \left(\int_{t_0}^T K_\mathcal{G}(s)ds \right)^{1/2} \leq \tilde{Z}B_2^{1/2}(T-t_0)^{1/3}.
\end{eqnarray*}
This shows that $H_\mathcal{G}(t)$ has a limit for $t \nearrow T$.
\end{proof}

\begin{proof}[Conclusion of the proof of (\ref{eq:thmJ-estm},\ref{eq:thmdotJ-estm}) in Theorem~\ref{thm:estm-coll-time}]
 From  Lemma~\ref{lem:HC-bnd}, it follows that all the sequence of lemmas which led to Lemma~\ref{lem:JdotJ-estm} is also satisfied for $J=J_\mathcal{G}$.
\end{proof}

\subsubsection{Proof of asymptotic estimates for $U_\mathcal{G}$ and $\ddot{J}_\mathcal{G}$ }
\label{sssec:sec-part}
Now, we aim to prove the following lemma.
\begin{lemma}
\label{lem:behavior of U}
    Under the same assumption as in Thm. \ref{thm:estm-coll-time},
    \begin{align*}
        \frac{U_\mathcal{G}(t)}{(T-t)^{-2/3}} \to \frac{2A}{9}, \\
        \frac{\ddot{J}_\mathcal{G}(t)}{(T-t)^{-2/3}} \to \frac{4A}{9}.
    \end{align*}
\end{lemma}
\begin{proof}
We primarily follow the approach in \cite{PS68}, switching to \cite{Wintner41} to apply a Tauberian theorem

From Lemmas~\ref{lem:ddj} and ~\ref{lem:HC-bnd}, it follows that
\begin{equation}
  \ddot{J}_\mathcal{G}=2U_\mathcal{G} + b(t), \label{eq:J''=2U+b}
\end{equation}
where $b(t)$ is a bounded continuous function.  Thus, to show that $\ddot{J}_\mathcal{G}$ is bounded, it is enough to establish the asymptotic bound for $U_\mathcal{G}$.

From already proved relation~\eqref{eq:thmdotJ-estm}, we have
\begin{equation}
   \frac{\dot{J}_{\mathcal{G}}(t)}{(T-t)^{1/3}} \to -\frac{4A}{3}, \quad t \nearrow T. \label{eq:dJG-asympt}
\end{equation}
Hence, $\dot{J}_\mathcal{G}(t) \to 0$ as $t \nearrow T$. Therefore, from \eqref{eq:J''=2U+b}, we have
\begin{equation*}
  -\dot{J}_\mathcal{G}(t)=2 \int_{t}^T U_\mathcal{G}(s)ds + O(T-t).
\end{equation*}
Let us denote
\begin{equation*}
  G(t)=-\int_{t}^T U_\mathcal{G}(s)ds,
\end{equation*}
so that
\begin{align*}
    \dot{J}_{\mathcal{G}}(t) = 2G(t)+ O(T-t).
\end{align*}
From Eq.~\eqref{eq:dJG-asympt}, it follows that
\begin{equation}
  \frac{G(t)}{(T-t)^{1/3}} \to  -\frac{2A}{3}, \quad t \nearrow T.  \label{eq:G-asympt}
\end{equation}
Obviously $\dot{G}(t)=U_\mathcal{G}(t)$. To obtain a bound on $\dot{G}(t)$, we need some bounds on $\ddot{G}(t)=\dot{U}_\mathcal{G}(t)$.

For this end we first prove Lemma~\ref{lem:dUU52} and then we will return to the proof of our lemma. 

\begin{lemma}
\label{lem:dUU52}
Under the assumptions of Theorem~\ref{thm:estm-coll-time}, there exists a constant $B$ such that
 \begin{equation}
   |\dot{U}_\mathcal{G}| \leq B U_\mathcal{G}^{5/2}.   \label{eq:U'G}
 \end{equation}
\end{lemma}
\begin{proof}
Let $m_0$ be the smallest mass in cluster $\mathcal{G}$ and let $r_{jk}=|q_k - q_j|$. We have
  \begin{eqnarray*}
    U_\mathcal{G} \geq \sum_{j <k, k,j \in \mathcal{G}} \frac{m_0^2}{r_{jk}}.
  \end{eqnarray*}
  Besides, we know that
  \begin{eqnarray*}
    -\dot{U}_{\mathcal{G}}=\sum_{j <k, k,j \in \mathcal{G}}\frac{m_k m_j}{r_{jk}^2} \dot{r}_{jk}
  \end{eqnarray*}
  Using the Cauchy-Schwarz inequality, we get
  \begin{eqnarray*}
   \left| \sum_{j <k, k,j \in \mathcal{G}}\frac{m_j m_k}{r_{jk}^2}   \dot{r}_{jk} \right| \leq \left( \sum_{j <k, k,j \in \mathcal{G}}\frac{m_j m_k}{r_{jk}^4} \right)^{1/2}
    \left( \sum_{j <k, k,j \in \mathcal{G}}m_j m_k \dot{r}^2_{jk}\right)^{1/2}.
  \end{eqnarray*}
  It is evident that there exist constants $C_1$ and $C_2$ depending on the masses such that
  \begin{eqnarray*}
     \sum_{j <k, k,j \in \mathcal{G}}\frac{m_j m_k}{r_{jk}^4} &\leq& C_1 U^4_\mathcal{G},\\
    \sum_{j <k, k,j \in \mathcal{G}}m_j m_k \dot{r}^2_{jk} &=&  \sum_{j <k, k,j \in \mathcal{G}}m_j m_k (\dot{q}_j - \dot{q}_k)^2 \leq C_2 K_\mathcal{G}.
  \end{eqnarray*}

Combining the above inequalities, we obtain
  \begin{eqnarray*}
    |\dot{U}_\mathcal{G}|^2 \leq C_1 C_2 U^4_\mathcal{G} K_\mathcal{G}=C_1C_2  U^4_\mathcal{G}(U_\mathcal{G} + H_\mathcal{G}).
  \end{eqnarray*}
From Lemma~\ref{lem:HC-bnd}, we know that $H_\mathcal{G}(t)$ is bounded. Since $U_{\mathcal{G}}(t) \to \infty$, we can find a constant $B$ such that
  \begin{equation*}
     |\dot{U}_\mathcal{G}|^2 \leq B U_\mathcal{G}^5.
  \end{equation*}

\end{proof}

Now we return to the proof of Lemma~\ref{lem:behavior of U}.

In preparation for using Tauberian Theorem~\ref{thm:tauberian}, following the approach of Wintner \cite{Wintner41}, we define
\begin{equation}
  F(\sigma)=G(t(\sigma))^3 ,\quad \text{ where} \quad \sigma = T-t.  \label{eq:F=G3}
\end{equation}

Rewriting Eq.~\eqref{eq:U'G} in terms of function $G$, we get
\begin{equation}
  |\ddot{G}| \leq B \dot{G}^{5/2}. \label{eq:ddG-estm}
\end{equation}
Denoting by $(\,)'$ the derivative with respect to $\sigma$, we have the following equalities.
\begin{eqnarray*}
  &F'=3G^2\dot{G} \frac{dt}{d\sigma} = -3G^2\dot{G}, \quad   F''=6G \dot{G}^2 + 3G^2 \ddot{G}, \\
   &G=F^{1/3}, \quad \dot{G}=\frac{-F'}{3 F^{2/3}}.
\end{eqnarray*}
Using Eq.~\eqref{eq:ddG-estm}, we obtain the following estimate for $|F''|$ in terms of $F'$ with some constant $C_3$.
\begin{eqnarray*}
  |F''| \leq 6 |F|^{1/3} \left(\frac{|F'|}{(3|F|^{2/3}}) \right)^2 + 3 |F|^{2/3} B \dot{G}^{5/2}=\frac{2 F'^2}{3|F|} + 3B |F|^{2/3}\left( \frac{|F'|}{3|F|^{2/3}}\right)^{5/2} \\
  = \frac{2 F'^2}{3 |F|} + \frac{3^{-3/2}B|F'|^{5/2}}{|F|}\leq C_3 \frac{F'^2 + |F'|^{5/2}}{|F|}.
\end{eqnarray*}
From  Eqs.~\eqref{eq:G-asympt} and \eqref{eq:F=G3},
\begin{equation*}
  \frac{F}{\sigma}\to -\left(\frac{2A}{3}\right)^3, \quad \sigma\to 0^+.
\end{equation*}
Therefore, there exists a constant $C$, such that,  for $\sigma$ close enough to $0$,
\begin{equation*}
 |F''|  \leq C \left(F'^2 + | F'| ^{5/2} \right)/\sigma.  \label{eq:19}
\end{equation*}
From Theorem~\ref{thm:tauberian}, we deduce that the following limit exists
\begin{equation*}
  \lim_{\sigma\to 0^+} F'(\sigma)=-\left(\frac{2A}{3}\right)^3 . \label{eq:limF'}
\end{equation*}

To complete the proof, we need to compute $\lim_{t\to T} U_\mathcal{G}(t)/(T-t)^{-2/3}$. By noticing that $t\nearrow T$ is equivalent to $\sigma \to 0^+$ and using Eq.~\eqref{eq:G-asympt}, we have
\begin{eqnarray*}
  \frac{U_\mathcal{G}(t)}{(T-t)^{-2/3}}=\frac{\dot{G}(t)}{(T-t)^{-2/3}}=\frac{-F'(\sigma(t))}{3 G^{2}(t)(T-t)^{-2/3}}=\frac{-1}{3}\frac{F'(\sigma(t))}{ \left(\frac{G(t)}{(T-t)^{1/3}}\right)^{2}}
  \to \frac{1}{3}\frac{\left(\frac{2A}{3}\right)^3}{\left(\frac{2A}{3}\right)^2}=\frac{2A}{9}.
\end{eqnarray*}
This completes the proof of Lemma~\ref{lem:behavior of U}.
\end{proof}

To finish the proof of Theorem~\ref{thm:estm-coll-time}, we only need to establish \eqref{eq:KG-bnds}.  This is a direct consequence of Eq.~\eqref{eq:U-bnds} and Lemma~\ref{lem:HC-bnd}.

\section{Convergence to the set of central configurations}
\label{sec:conv-cc}

Our goal is to prove  the shape of the configuration of bodies going toward collision approaches the set of central configurations. For total collisions, this was first proved by Chazy \cite{Ch18}.
For partial collisions such statement appeared in a paper by Saari \cite{Saari84}, but we have doubts about the use of Tauberian theorem in its proof.

Therefore, we provide a proof based on the adaption of Wintner's proof of this fact for total collision \cite[\S 361-365]{Wintner41}
\begin{theorem}
\label{thm:conver-to-cc}
The same assumptions as in Theorem~\ref{thm:estm-coll-time}.

Then $\left(\frac{q_i-L_\mathcal{G}}{|T-t|^{2/3}}\right)_{i \in \mathcal{G}}$  approaches the set of central configurations.
\end{theorem}
\begin{proof}
To follow the arguments from Wintner's book \cite{Wintner41} we will assume that collision happens for $t \to 0^+$, i.e. we go backward in time. 

 We  assume that $L_\mathcal{G}=0$ (i.e. partial collision happens at $0$).
We introduce rescaled coordinates
\begin{equation*}
  \xi_i(t)=q_i(t)/t^{2/3}, \quad i \in \mathcal{G}.
\end{equation*}
We want to prove that vector $\xi(t)$ approaches the set of central configurations.

From  equations (\ref{eq:second-order-formulation})  for $q_i$ in our cluster $\mathcal{G}$ we obtain the following equations for $\xi_i$
\begin{eqnarray*}
  m_i\left( t^{2/3} \ddot{\xi}_i + \frac{4}{3}t^{-1/3} \dot{\xi}_i - \frac{2}{9}t^{-4/3}\xi_i \right)=t^{-4/3}\sum_{j \in \mathcal{G}, j \neq i} \frac{m_i m_j (\xi_j - \xi_i)}{|\xi_i - \xi_j|^3} + \frac{\partial U_{ext}}{\partial q_i},
\end{eqnarray*}
which after multiplication by $t^{4/3}$ becomes  (with $U(\xi)=U_\mathcal{G}(\xi)$)
\begin{eqnarray}
   m_i\left( t^{2} \ddot{\xi}_i + \frac{4}{3}t \dot{\xi}_i - \frac{2}{9}\xi_i \right)=\sum_{j \in \mathcal{G}, j \neq i} \frac{m_i m_j (\xi_j - \xi_i)}{|\xi_i - \xi_j|^3} + t^{4/3} f(t) = \frac{\partial U}{\partial \xi_i}(\xi) +  t^{4/3} f(t), \label{eq:resc-N-body} \\
   f(t)=O(1), \quad \dot{f}(t)=O(t^{-1/3}).  \label{eq:resc-N-body-bndsg}
\end{eqnarray}
Bounds in (\ref{eq:resc-N-body-bndsg}) are a consequence of the fact that $U_{\text{ext}}$ and its partial derivatives are bounded when we are going to the partial collision.
This gives $f(t)=O(1)$.
To obtain bound on  $\dot{f}(t)$ note that
\begin{equation*}
\dot{f}(t)=\sum_{j=1}^n \frac{\partial }{\partial q_j}\left(\frac{\partial U_{ext}}{\partial q_i}\right)(q(t))\dot{q}_j(t)
\end{equation*}
and then estimate follows from bound $\dot{q}_j(t)=O(t^{-1/3})$ (see Eq. (\ref{eq:KG-bnds}) in Theorem~\ref{thm:estm-coll-time} for all colliding clusters, while for non-colliding bodies $\dot{q}_j(t)=O(1)$).

Now we introduce a new time variable $\tau=-\ln t$ (i.e. $e^{-\tau}=t$, $\frac{d t}{d\tau}=-t$ and $\tau \to \infty$ when $t \to 0$). We will denote derivatives with respect to $\tau$ by dash $'$.

We have for any function $x$
\begin{eqnarray}
  \dot{x}=- x'/t, \quad \ddot{x}=x''/t^2 + x'/t^2,  \label{eq:diff-rel}
\end{eqnarray}
therefore (\ref{eq:resc-N-body}) becomes
\begin{eqnarray}
   m_i \left(\xi_i'' - \frac{\xi_i'}{3} - \frac{2 \xi_i}{9}\right)= \frac{\partial U}{\partial \xi_i}(\xi) +  e^{-4 \tau /3} f(e^{-\tau}), \quad i \in \mathcal{G} \label{eq:19-1}
\end{eqnarray}
Our goal is to show that
\begin{equation}
\xi'' \to 0, \quad  \xi' \to 0 \qquad \mbox{for } \quad \tau \to \infty,  \label{eq:xder-lim}
\end{equation}
 so that, in the limit, we obtain the equation
\begin{equation}
 -\frac{2 m_i \xi_i}{9}= \frac{\partial U}{\partial \xi_i}(\xi), \label{eq:CC-eq-xi}
\end{equation}
which is the equation for central configurations.

To achieve (\ref{eq:xder-lim}), it is enough to show that $\xi' \to 0$ and $\xi'''$ is bounded, then using  Tauberian theorem~\ref{thm:taubwin} from Appendix~\ref{sec:TaubThm}, we will infer that also $\xi'' \to 0$.

Now we will work towards estimating $\xi'$ and $\xi'''$. This will be achieved by control of $J_\mathcal{G}$ and its derivatives obtained in Theorem~\ref{thm:estm-coll-time}, equation (\ref{eq:resc-N-body}) and equation for $\ddot{J}$ obtained in Lemma~\ref{lem:ddj}.

Since  the collision happens for $t \to 0^+$  the results of Theorem~\ref{thm:estm-coll-time} can be rewritten
as follows
\begin{align}
  \frac{J_{\mathcal{G}}(t)}{t^{4/3}} \to A, \quad t \to 0^+  \label{eq:thmJ-estm-r} \\
  \frac{\dot{J}_{\mathcal{G}}(t)}{t^{1/3}} \to \frac{4A}{3}, \quad  t \to 0^+  \label{eq:thmdotJ-estm-r} \\
    \frac{\ddot{J}_\mathcal{G}(t)}{t^{-2/3}} \to \frac{4}{9}  A,  \quad  t \to 0^+.   \label{eq:J''-bnds}
\end{align}

From now on, we will drop subscript in  $J_{\mathcal{G}}$ and just write $J$, similarly $U_\mathcal{G} \to U$ and other bodies are acting on bodies in cluster $\mathcal{G}$
through $U_{\text{rel}}$. Since $L_\mathcal{G}=0$, so then
\begin{equation*}
  J=\sum_{i \in \mathcal{G}} m_i q_i^2.  \label{eq:def-JG}
\end{equation*}

Easy computations show that (we use (\ref{eq:thmJ-estm-r},\ref{eq:thmdotJ-estm-r},\ref{eq:J''-bnds}))
\begin{eqnarray}
  t\frac{d}{dt} (t^{-4/3}J) \to 0 , \quad t \to 0^+,  \label{eq:14-2} \\
  t^2\frac{d^2}{d t^2} (t^{-4/3}J)\to 0 , \quad t \to 0^+. \label{eq:14-3}
\end{eqnarray}

From Lemmas~\ref{lem:ddj} and~\ref{lem:HC-bnd} we know that
\begin{equation}
  \ddot{J}(t)=2U(t) + g(t),   \label{eq:J''=2U+g}
\end{equation}
where $g(t)$ is bounded function.

Now we will  rewrite equation (\ref{eq:J''=2U+g}) in terms of new variables.
Since   
\begin{equation}
J(q(t))=t^{4/3}J(\xi(t)). \label{eq:Jxi}
\end{equation}
and
\begin{eqnarray*}
  \frac{d^2}{dt^2} \left( t^{4/3}J(\xi) \right)=\frac{4}{9}t^{-2/3}J(\xi) + \frac{8}{3}t^{1/3}\dot{J}(\xi) +t^{4/3}\ddot{J}(\xi).
\end{eqnarray*}
Using relations (\ref{eq:diff-rel})  we obtain
\begin{eqnarray}
 \ddot{J}=\frac{d^2}{dt^2} \left( t^{4/3}J(\xi) \right)= \frac{4}{9}t^{-2/3}J(\xi) - \frac{8}{3}t^{-2/3}J'(\xi) +t^{-2/3}(J'(\xi) + J''(\xi))  \nonumber \\
 =t^{-2/3}\left( \frac{4}{9}J(\xi)  - \frac{5}{3}J'(\xi)+ J''(\xi)\right),  \label{eq:ddJ-newvar}
\end{eqnarray}
hence equation (\ref{eq:J''=2U+g}) becomes
\begin{equation}
   J''(\xi)   - \frac{5}{3}J'(\xi)+  \frac{4}{9}J(\xi)= 2 U(\xi) + e^{-2\tau/3} g(e^{-\tau}).  \label{eq:20-2}
\end{equation}
From (\ref{eq:thmJ-estm-r}) and (\ref{eq:Jxi}) it follows immediately that
\begin{equation}
  J(\xi(\tau)) \to A, \quad \tau \to \infty.  \label{eq:21-1}
\end{equation}

Since (we use (\ref{eq:diff-rel}))
\begin{eqnarray*}
  J'(\xi(\tau))= -t \dot{J}(\xi(t))=-t\frac{d}{dt}\left(t^{-4/3}J(q(t)) \right)
\end{eqnarray*}
then by (\ref{eq:14-2})  it follows that
\begin{equation}
  J'(\xi) \to 0, \quad \tau \to \infty.  \label{eq:21-2}
\end{equation}

Since (we use (\ref{eq:diff-rel}))
\begin{eqnarray*}
  J''(\xi(\tau))=t^2 \ddot{J}(\xi(t)) + J'(\xi(\tau))=t^2 \frac{d^2}{dt^2}\left( t^{-4/3}J(q(t)) \right) + J'(\xi(\tau))
\end{eqnarray*}
then from (\ref{eq:14-3}) and (\ref{eq:21-2}) we obtain
\begin{equation}
  J''(\xi) \to 0, \quad \tau \to \infty.  \label{eq:21-3}
\end{equation}

Passing to the limit $\tau \to \infty$ in (\ref{eq:20-2}) and using (\ref{eq:21-1},\ref{eq:21-2},\ref{eq:21-3}) we obtain
that
\begin{equation}
   \frac{2}{9}J(\xi) - U(\xi) \to 0, \quad \tau \to \infty \quad (t \to 0^+). \label{eq:17-1}
\end{equation}

Moreover, since $J(\xi)$ has a finite limit, $U(\xi)$ is bounded when going to the collision. This implies that  there exists a constant $C>0$
such that
\begin{equation}
  |\xi_k - \xi_j| \geq C, \quad i\neq j, \quad \forall t.  \label{eq:xi-dist-bnd}
\end{equation}

From Lemma~\ref{lem:HC-bnd} we have
\begin{equation}
  \frac{1}{2}\sum_{i \in \mathcal{G}} m_i \dot{q}_i^2 - U(q) = h(t) \label{eq:ener-eq-pr}
\end{equation}
where $h(t)$ is some bounded function.

Since
\begin{equation*}
  \dot{q}_i=-t^{-1/3}\left(\xi_i' - \frac{2}{3}\xi_i \right)
\end{equation*}
then  relation (\ref{eq:ener-eq-pr})  becomes
\begin{eqnarray}
  \frac{1}{2}\sum_{i \in \mathcal{G}} m_i \left(\xi_i' - \frac{2}{3}\xi_i \right)^2 - U(\xi)=e^{-2\tau/3}h(e^{-\tau}).  \label{eq:20-1}
\end{eqnarray}

Next, we show that $\xi' \to 0$. Relation (\ref{eq:20-1}) gives
\begin{eqnarray*}
  \frac{1}{2}\sum_i m_i \xi'^2 - \frac{1}{6}J'(\xi) + \frac{2}{9}J(\xi) - U(\xi)=e^{-2\tau/3}h(e^{-\tau})
\end{eqnarray*}
Now passing to the limit $\tau \to \infty$ and using (\ref{eq:17-1},\ref{eq:21-2}) we obtain
\begin{equation}
  \xi' \to 0, \quad \tau \to \infty. \label{eq:22-1}
\end{equation}

It is easy to see that $\xi$, $\xi'$ and $\xi''$ are bounded. In case of $\xi''$ this follows from (\ref{eq:19-1}) , the boundedness of $\frac{\partial U}{\partial \xi_i}$
is granted by (\ref{eq:xi-dist-bnd}).

Now, we are ready to show that $\xi'''$ is bounded.
For this end, we differentiate equation (\ref{eq:19-1}).  Observe that from (\ref{eq:resc-N-body-bndsg}) we obtain
for some constant $C$
\begin{eqnarray*}
\left(e^{-4 \tau /3} f(e^{-\tau})\right)' = \frac{-4}{3}e^{-4 \tau /3} f(e^{-\tau}) - e^{-7\tau/3} \dot{f}(e^{-\tau}) , \\
\left|\left(e^{-4 \tau /3} f(e^{-\tau})\right)'\right| \leq   e^{-4 \tau /3} C + e^{-7\tau/3} C e^{\tau/3} = 2C   e^{-4 \tau /3}.
\end{eqnarray*}
Due to (\ref{eq:xi-dist-bnd}) we have that partial derivatives of $U(\xi)$ are bounded, hence from (\ref{eq:22-1}) we obtain that $\left(U(\xi(\tau))\right)'$ are bounded.
Since $\xi''$ and $\xi'$ are bounded, hence $\xi'''$ is also bounded.

Now, the application of the  Tauberian Theorem~\ref{thm:taubwin}
 gives that $\xi'' \to 0$.

Passing to the limit $\tau \to \infty$ in equation (\ref{eq:19-1}) gives (\ref{eq:CC-eq-xi}). This finishes the proof.
\end{proof} 

\section{Two Tauberian theorems}
\label{sec:TaubThm}
The  following Tauberian theorem was stated and proved in \cite{Wintner41} (in \S 338 and 338bis), a particular case of this theorem with $\nu_0=0$ appeared as Theorem 1 in \cite{B39}.
 \com{PZ: we did not like the proof there }
\begin{theorem}
\label{thm:tauberian}
 Let $F:(t_0,T] \to \mathbb{R}_+$  is $C^2$ and for some continuous increasing function $g:\mathbb{R}_+ \to \mathbb{R}_+$
 holds  that
 \begin{equation*}
   |F''(t)| \leq g(|F'(t)|)/(t-t_0), \quad t \in (t_0,T]
 \end{equation*}
 \begin{eqnarray*}
   F \sim \nu_0(t-t_0),
 \end{eqnarray*}
 with some $\nu_0 \in \mathbb{R}$,
 then
 \begin{equation*}
    \lim_{t \to t_0} F'(t) = \nu_0.
 \end{equation*}
\end{theorem}

We use also the following Tauberian Theorem
 \begin{theorem}\cite[\S 362]{Wintner41}
\label{thm:taubwin}
Assume $g:(0,\infty) \to \mathbb{R}$ is $C^2$, $\lim_{u \to \infty} g(u)$ exists and $|\ddot{g}(u)| \leq M$ for all $u \in (0,\infty)$. Then
\begin{equation*}
  \lim_{u \to \infty} \dot{g}(u)=0.
\end{equation*}
 \end{theorem}



\section*{Conflicts of Interest Statement}
The authors declare that they have no conflict of interest.

\section*{Data availability statement}
No datasets were generated or analysed during the current study.

\section*{Acknowledgement}
RGS was supported by the Priority Research Area SciMat under the program Excellence Initiative - Research University at the Jagiellonian University in Krak\'ow.

PZ was  supported by Polish NCN grants 2019/35/B/ST1/00655  and 2021/41/B/ST1/00407.

AG was supported by Polish NCN grant 2021/41/B/ST1/00407.
\bibliography{references}

\begin{thebibliography}{McG78}

\bibitem[AK12]{AlbouyKaloshin12}
A.~Albouy and V.~Kaloshin.
\newblock Finiteness of central configurations of five bodies in the plane.
\newblock \emph{Annals of Mathematics}, 176:535--588, 2012.

\bibitem[Boa39]{B39}
R.~P. Boas, Jr.
\newblock A {T}auberian theorem connected with the problem of three bodies.
\newblock \emph{Amer. J. Math.}, 61:161--164, 1939.

\bibitem[Cha18]{Ch18}
J.~Chazy.
\newblock Sur certaines trajectoires du probleme des $n$ corps.
\newblock \emph{Bull. Astron.}, 35:321--389, 1918.

\bibitem[CZ15]{CZ15}
M.~J. Capi\'nski and P.~Zgliczy\'nski.
\newblock Geometric proof for normally hyperbolic invariant manifolds.
\newblock \emph{Journal of Differential Equations}, 259(11):6215--6286, 2015.

\bibitem[ElB90]{E}
M.~ElBialy.
\newblock Collision singularities in the n-body problem.
\newblock \emph{SIAM J. Math. Anal.}, 21(6):1563--1593, 1990.

\bibitem[HM06]{HamtonMoeckel06}
M.~Hampton and R.~Moeckel.
\newblock Finiteness of relative equilibria of the four-body problem.
\newblock \emph{Inventiones Mathematicae}, 163:289--312, 2006.

\bibitem[McG74]{McGehee74}
R.~McGehee.
\newblock Triple collision in the collinear three-body problem.
\newblock \emph{Invent Math}, 27:191--227, 1974.

\bibitem[McG78]{McGehee78}
R.~McGehee.
\newblock Singularities in classical celestial mechanics.
\newblock \emph{Proc. of Int. Cong. Math., Helsinki}, pages 827--834, 1978.

\bibitem[MM24]{Moeckel2023}
R.~Moeckel and R.~Montgomery.
\newblock No infinite spin for planar total collision.
\newblock \emph{Journal of the American Mathematical Society}, 01 2024.

\bibitem[Moe23]{MoeckelTCSC}
R.~Moeckel.
\newblock Total collision with slow convergence to a degenerate central
  configuration.
\newblock \emph{Regular and Chaotic Dynamics}, 28:533--542, 2023.

\bibitem[MZ19]{MZ19}
M.~Moczurad and P.~Zgliczy\'nski.
\newblock Central configurations in planar $n$-body problem for $n=5,6,7$ with
  equal masses.
\newblock \emph{Celestial Mechanics and Dynamical Astronomy}, 131:46, 2019.

\bibitem[Pal76]{Palmore76}
J.~Palmore.
\newblock Measure of degenerate relative equilibria i.
\newblock \emph{Annals of Mathematics}, 104:421--431, 1976.

\bibitem[PS68]{PS68}
H.~Pollard and D.~G. Saari.
\newblock Singularities of the $n$-body problem {I}.
\newblock \emph{Arch. Rational Mech. Anal.}, 30:263--269, 1968.

\bibitem[Saa71]{Saari71}
D.~G. Saari.
\newblock Expanding gravitational systems.
\newblock \emph{TAMS}, 156:219--214, 1971.

\bibitem[Saa84]{Saari84}
D.~G. Saari.
\newblock The manifold structure for collision and for hyperbolic-parabolic
  orbits in the $n$-body problem.
\newblock \emph{Journal of Differential Equations}, 55(3):300--329, 1984.

\bibitem[SH81]{SH}
D.~Saari and N.~Hulkower.
\newblock On the manifolds of total collapse orbits and of completely parabolic
  orbits for the $n$-body problem.
\newblock \emph{Journal of Differential Equations}, 41:27--43, 1981.

\bibitem[Shu70]{Sh70}
M.~Shub.
\newblock Diagonals and relative equilibria, {A}ppendix to {S}male's paper.
\newblock \emph{Springer Lecture Notes in Math.}, 197:199--201, 1970.

\bibitem[Sma98]{Smale98}
S.~Smale.
\newblock Mathematical problems for the next century.
\newblock \emph{The Mathematical Intelligencer}, 20:7--15, 1998.

\bibitem[Spe70]{Sp}
H.~Sperling.
\newblock On the real singularities of the n-body problem.
\newblock \emph{J. Reine Angew. Math.}, 245:15--40, 1970.

\bibitem[Srz94]{Srz}
R.~Srzednicki.
\newblock Periodic and bounded solutions in blocks for time-periodic
  nonautonomuous ordinary differential equations.
\newblock \emph{Nonlin. Analysis, TMA}, 22:707--737, 1994.

\bibitem[Sun13]{S13}
K.~F. Sundman.
\newblock M\'emoire sure le probl\'eme des trois corps.
\newblock \emph{Acta Math.}, 36(1):105--179, 1913.

\bibitem[Wa{\.z}47]{waz}
T.~Wa{\.z}ewski.
\newblock Sur un principe topologique de l'examen de l'allure asymptotique des
  int\'egrales des \'equations diff\'erentielles ordinaires.
\newblock \emph{Ann. Soc. Polon. Math.}, 20:279--313, 1947.

\bibitem[Wig94]{Wiggins}
S.~Wiggins.
\newblock \emph{Normally Hyperolic Invariant Manifolds in Dynamical Systems},
  volume 105.
\newblock 01 1994.
\newblock ISBN 978-1-4612-8734-6.

\bibitem[Win41]{Wintner41}
A.~Wintner.
\newblock \emph{The Analytical Foundations of Celestial Mechanics}.
\newblock Princeton Univ. Press, 1941.

\bibitem[Zgl17]{Zsegments}
P.~Zgliczy\'nski.
\newblock Topological shadowing and the {G}robman-{H}artman theorem.
\newblock \emph{Topological Methods in Nonlinear Analysis}, 49:1, 11 2017.

\end{thebibliography}
\end{document}